\documentclass[11pt]{article}
\usepackage[a4paper,margin=3cm]{geometry}
\usepackage{amsmath,amssymb,amsthm,mathtools}
\usepackage{hyperref}
\usepackage[english]{babel}
\usepackage[nameinlink,capitalize]{cleveref}
\usepackage{csquotes}
\usepackage[english]{babel}
\usepackage{enumitem}
\usepackage{shuffle}
\usepackage{verbatim}
\usepackage[blocks]{authblk}
\usepackage[
  backend=biber,
  style=alphabetic,
  maxbibnames=99
]{biblatex}
\crefname{appendix}{Appendix}{Appendices}
\Crefname{appendix}{Appendix}{Appendices}

\theoremstyle{definition}
\newtheorem{definition}{Definition}
\newtheorem{lemma}{Lemma}
\newtheorem{remark}{Remark}
\newtheorem{claim}{Claim}
\newtheorem{theorem}{Theorem}

\newtheorem{example}{Example}

\newtheorem{corollary}{Corollary}
\newtheorem{proposition}{Proposition}

\title{Recovering the Length of a Path from its Signature}
\date{\today}
\author{Felix Medwed\thanks{%
    Institute of Mathematics, University of Potsdam, D-14469 Potsdam, Germany\\
    E-mail: \href{mailto:felix.medwed@uni-potsdam.de}
                 {\texttt{felix.medwed@uni-potsdam.de}}}}

\begin{document}

\maketitle

\begin{abstract}
    We prove the length conjecture of J.~Chang, T.~Lyons and H.~Ni for paths of bounded variation for arbitrary real Banach spaces and any reasonable tensor (algebra) norms. More precisely, we establish the missing lower bound and in this way show that the length of a tree-reduced path is recovered from the asymptotic behaviour of its normalised signature. The main technical argument uses the norm duality and shuffle algebra to express natural powers of polygonal length approximations as combinations of linear functionals acting on restricted signatures, which are then glued by suitable integrations into a linear functional of controlled growth acting on the full signature. In particular, the full proof combines tools of classical analysis and functional analysis, such as localisation by continuous cutoff functions, Stone--Weierstrass approximation of continuous functions by signature polynomials, the Hahn-Banach theorem and a support detection lemma for compactly supported real signed Borel measures.
\end{abstract}

\tableofcontents

\section*{Preamble}

The central arguments in the proof of the length conjecture were generated with the assistance of AI (ChatGPT 5.6 Sol) back in summer of 2026. The starting point was the appendix of the article \cite{kerger2026closingoraclecomplexitygapderivativefree}, following \cite{OpenAIcdn}, and which the author believed could be adapted to the length conjecture. The author then asked the AI to modify the corresponding prompt accordingly, so that it could serve as a potential research prompt for the length conjecture. The arguments produced by the AI were then subsequently turned by the author into a concise proof, with the necessary restructering, verifications, corrections and generalisations. Since then, the article has undergone numerous revisions. The proof may therefore be viewed as the result of an intricate intertwining of human and AI contributions. Furthermore, the author presents this proof in the belief that the employed techniques may prove useful beyond the present exposition. In that sense, a further investigation using other techniques to show the length conjecture remains worthwhile.

On request the author can provide the raw research prompt and raw ChatGPT proof.

\section{Motivation and historical context}
\label{sec:motivation-and-historical-context}
The \emph{signature transformation}, or simply \emph{signature}, of a continuous bounded variation path \(\gamma \colon [a,b] \to V\) for \(a\neq b\) is the formal tensor series formed by its iterated integrals.
Such iterated integral transformations were considered for finite dimensional \(V\) and piecewise continuously differentiable paths by K.T.~Chen \cite{Chen1954} to construct a cohomology theory of loop spaces. Furthermore, similar transformations emerged in the work of Dyson under the name time-ordered expansions in perturbative quantum field theory. It has since sparked the \emph{Theory of Rough Paths}, introduced by T.~Lyons \cite{Lyons1998}, which inspired M.~Hairer's theory of regularity structures.\\
In their seminal work \cite{HamblyLyons2010}, B.~Hambly and T.~Lyons proved, for \(V\) being a finite-dimensional vector space and for any continuous bounded variation path \(\gamma\), that the signature faithfully encodes the geometric data of a path up to \emph{tree-like equivalence}. The follow-up work \cite{BoedihardjoGengLyonsYang2016} of H.~Boedihardjo, X.~Geng, T.~Lyons and D.~Yang then extended this result to arbitrary real Banach spaces and a rougher class of paths. Informally, the equivalence accounts for excursions in which the path returns along itself, leading to a full cancellation in its signature. Furthermore, up to the fixation of a starting point and parametrisation of the considered path \(\gamma\), there exists a unique representative that exhibits no retracing, referred to as the \emph{tree-reduced} representative \(\gamma_{\mathrm{red}}\) associated to \(\gamma\).\\
In the same work, B.~Hambly and T.~Lyons established, for suitable choices of tensor norms, that the \emph{length} (or \emph{variation}) \(L(\gamma)\) of a continuously differentiable path, parametrised to arc length, can be recovered from its signature through the asymptotic formula
\begin{equation} \label{eq:asymptotics-signature}
    \lim_{n \to \infty} \left( n! \left\| \int_{a < t_1 < ... < t_n < b} \mathrm{d}\gamma(t_1) \otimes ... \otimes \mathrm{d}\gamma(t_n) \right\|_n \right)^{1/n} = L(\gamma).
\end{equation}
It was then explicitely conjectured by J.~Chang, T.~Lyons and H.~Ni \cite{ChangLyonsNi2018} with a follow-up \cite{BoedihardjoGeng2019} by H.~Boedihardjo and X.~Geng, that this recovery persists also without imposing additional regularity, given the choice of \emph{reasonable tensor (algebra) norms} and assuming the path is tree-reduced. This is now referred to as the \emph{length conjecture}, where the corresponding limit \eqref{eq:asymptotics-signature} was shown to exists, ref.~\cite[Theorem 3 with corrigendum]{ChangLyonsNi2018} and \cite[Corollary 1]{BoedihardjoGeng2019}.

The conjecture was subsequently proved by H.~Boedihardjo and X.~Geng \cite{BoedihardjoGeng2023} for a class of planar paths under suitable control over singularities. They suggested a full proof to \textquote{[...] capture of the degree of tree-reducedness at a quantitative level, which is at the same time convenient to be utilised in the signature analysis}. 
Such an approach faces the difficulty that it must relate the absence of tree-like pieces to estimates involving the variation measure of the path and functionals of its signature, while allowing arbitrarily fine oscillations without imposing additional regularity. A further challenge arises in passing from \(\mathbb R^d\) to a general real Banach space, where coordinate maps must be replaced by suitable continuous linear functionals, whose choice must remain compatible with the lower bounds required in the signature analysis.

The proof presented here neither requires a quantitative description of tree reducedness of the kind considered in \cite{BoedihardjoGeng2023}, nor does it impose a finite-dimensional restriction on the underlying Banach space. We therefore bypass the explicit analysis of the local behaviour of the path. The main observations underlying the proof are the dual formulation of the length of polygonal approximations, a localisation argument that lets us retain the lower bound, the relation of natural powers of the dual to restricted signatures, a support detection lemma and lastly gluing the restricted signatures back to the full signature together with the identification of a suitable family of linear functionals with controlled growth. In particular, the key approach towards the lower bound is, in spirit, of the form considered in \cite{BoedihardjoGeng2023}, but does not rely on a fixed path development, e.g.~the \(\mathrm{SL}_2\)-development in \cite{BoedihardjoGeng2023}.

In addition, the identification step for the suitable family of linear functionals admits an algebraic formulation in terms of half-shuffles \cite{SalviDiehlLyonsPreiss2023} and commutative Rota-Baxter algebras of zero weight \cite{Kurusch2007,GaoGuoZhang2021}. It is therefore an observation, that the asymptotic control of the suitable family of linear functionals is induced by the underlying algebraic structure. However, we will refrain from introducing this perspective throughout for accessibility and choose to write out the involved expressions more explicitly.

Finally, we emphasise that the individual arguments employed in the proof are largely classical in rough path theory, analysis, functional analysis and probability theory. Nevertheless, to the best of our knowledge, their particular combination in the present setting has not yet been exploited. In addition, the proof uses a \textquote{Paley--Wiener type perspective}, already present in previous work on signatures such as \cite{LyonsSidorova2006,LyonsXu2017,BoedihardjoGeng2019Brownian,BoedihardjoGengSouris2020}, which we expect will guide even further proofs that concern signatures.
A further strength of the proof lies in its potential straight forward adaptation to higher signatures, such as surface signatures \cite{lee2024surfacesignature} through suitable modifications of the localisation and approximation arguments, as well as other choices of Hopf algebras. At last, we remark that the guiding principle of gluing local information into a global object underlying this proof is a well-established method in many areas of mathematics.

\section{Preliminaries}
\label{sec:introduction}

We denote by \(\mathbb{N}_0 = \mathbb{N} \cup \{0\}\) the natural numbers with zero. We assume throughout the real numbers \(\mathbb{R}\) to be equipped with their standard topology induced by the absolute value \(|\cdot|\) and the corresponding Borel \(\sigma\)-algebra. Let \(J:=[a,b]\) be a closed interval of the real numbers and \(m \in \mathbb{N}\). Then, we define
\[
    \Delta_{m}(J) := \{ (s_1,...,s_{m}) \mid a < s_1 < ... < s_{m} < b \} \subset J^{m}
\]
to be the \emph{\(m\)th standard simplex over \(J\)}.

For any topological spaces \(X\) and \(Y\) we consider, unless stated otherwise, the product space \(X \times Y\) to be equipped with the respective product topology and any subset \(S \subseteq X \times Y\) with the corresponding subspace topology.
Then, recall that the closure \(\overline{\Delta_{m}}\) is a compact set of \(J^{m}\). \\

We refer the reader for the subsequent preliminaries to the vast literature on rough path theory, such as \cite{Lyons2007}, \cite{Friz_Victoir_2010_fullbook} and \cite{FrizHairer2020}. For a more geometric take, see also \cite{Schmeding_2022}.


\subsection{Tensor series algebra and reasonable tensor norms}

Throughout, fix a real Banach space \(V\) with norm \(\|\cdot\|\). We will understand \(V\) as a metric space with induced metric \(\|\cdot- \cdot\|\) or topological vector space with the induced open ball topology when appropiate. 
Let, for all \(n \in \mathbb{N}_0\), \(V^{\otimes_a n}\) be the algebraic tensor powers of \(V\) with the convention that \(V^{\otimes_a 0} := \mathbb{R}\) and \(V^{\otimes_a 1}\) is identified with \(V\). Since each algebraic tensor power is naturally a real vector space, we may consider each algebraic tensor power to be equipped with a norm \(\|\cdot\|_n\) where by convention \(\|\cdot\|_0 := |\cdot|\) and \(\|\cdot\|_1 := \|\cdot\|\). Set \(V^{\otimes n} := \overline{V^{\otimes_a n}}^{\|\cdot\|_n}\) to be the completed tensor powers. 
Then, let
\[
    T((V)) = \prod_{n = 0}^\infty V^{\otimes n}
\]
be the formal tensor series algebra under the formal tensor multiplication \(\otimes\), which satisfies the two-sided distributive property
\[
    (u+v) \otimes w = u \otimes v + v \otimes w, \qquad u \otimes (v + w) = u \otimes v + u \otimes w.
\]
Let \(\pi_n \colon T((V)) \to V^{\otimes n}\) be the projection of a tensor series to its \(n\)th component and \(\pi_{\leq n} \colon T(V) \to T^n(V) := \oplus_{k = 0}^n V^{\otimes k}\) be the \(n\)th level truncation of \(T((V))\). For an element \(x \in T((V))\) we will directly index \(x_n := \pi_n(x)\) and \(x_{\leq n} := (x_0,x_1,...,x_n) \in T^n(V)\) as well as identify \(V^{\otimes n} \subset T^n(V) \subset T((V))\) under its canonical inclusion throughout the context. We then call any \(V^{\otimes n}\), understood under the canonical inclusion, a \emph{homogeneous component}. 
One may instead of the formal tensor series, consider the completed tensor algebra with respect to the Frechét topology induced by the semi-norms \((\|\pi_n(\cdot)\|_n)_{n \in \mathbb{N}_0}\). However, we will not utilise this explicitely.
Furthermore, we consider the \emph{continuous dual spaces} \(V^*\), \((V^{\otimes 2})^*\), ... to \(V\), \(V^{\otimes 2}\), ... with the their induced (operator) norms 
\[
    \|\phi\|_{n,*} := \sup_{x \in V^{\otimes n} \backslash \{0\}} \frac{|\phi_n(x)|}{\|x\|_n} = \sup_{\|x\|_n \leq 1} |\phi_n(x)|
\]
respectively and will employ the dual pairing notation \(\phi(x) \equiv \langle \phi, x \rangle\) for linear functionals \(\phi\) when convenient.
Recall that a linear functional is continuous if and only if it is bounded. 
The projections \(\pi_n\) onto the homogeneous components and the truncated tensor algebras allow us to identify their duals canonically with subspaces of the dual of the tensor series algebra.
In anticipation of the result, we consider each tensor level with the choice of \emph{reasonable tensor (algebra) norms} defined as follows.
\begin{definition}[Reasonable tensor (algebra) norm, comp.~\protect{\cite{ChangLyonsNi2018}}] \label{def:reasonable-tensor-norms}
    A family of norms \((\|\cdot\|_n)_{n \in \mathbb{N}}\) for the normed vector spaces \(V,V^{\otimes 2}, ... , V^{\otimes n},...\) respectively is called \emph{reasonable} if
    \begin{enumerate}[label = \arabic*)]
        \item \label{def:reasonable-tensor-norms:banach} the \emph{Banach inequality}
        \[
            \|v \otimes w\|_{n+m} \leq \|v\|_n \|w\|_m, \qquad \forall v\in V^{\otimes n}, \forall w \in V^{\otimes m}, \forall (n,m) \in \mathbb{N}^2,
        \]
        \item \label{def:reasonable-tensor-norms:permutation} for any \(n \in \mathbb{N}\), any permutation \(\sigma \in \mathrm{Sym}(n)\) and tensor \(u \in V^{\otimes n}\)
        \[
            \|u \cdot \sigma\|_n = \|u\|_n,
        \]
        where \(\_\cdot \sigma\) denotes the right action of \(\sigma \in \mathrm{Sym}(n)\) on \(V^{\otimes n}\) obtained by linearly extending
        \[
            (v_1\otimes\cdots\otimes v_n)\cdot \sigma := v_{\sigma(1)}\otimes\cdots\otimes v_{\sigma(n)},
        \]
        and
        \item \label{def:reasonable-tensor-norms:dual-banach} the \emph{dual Banach inequality}
        \[
            \|\phi \otimes \psi\|_{n+m,*} \leq \|\phi\|_{n,*} \|\psi\|_{m,*}, \qquad \forall \phi\in \left(V^{\otimes n}\right)^*, \forall \psi \in \left(V^{\otimes m}\right)^*, \forall (n,m) \in \mathbb{N}^2
        \]
    \end{enumerate}
    hold.
\end{definition}
\begin{remark}
    We emphasise that all of these conditions are used in the provided proof of the length conjecture. In particular, the Banach inequality \ref{def:reasonable-tensor-norms:banach} is invoked in the derivation of the upper bound \eqref{eq:Universal-Upper-Bound}, while the symmetry \ref{def:reasonable-tensor-norms:permutation} and the dual Banach inequality \ref{def:reasonable-tensor-norms:dual-banach} enter in the proof of \cref{lem:Growth-bound-of-functionals}, which provides the essential control on the growth of a suitable family of functionals.
\end{remark}
\begin{remark}
    It is an exercise using the Hahn--Banach theorem \cite[§23.18, (HB7), p.~617]{Schechter1997}, see also \cref{lem:norm-functional-association} later on, that reasonable tensor norms are cross-norms \cite{ChangLyonsNi2018}. 
    Hence, \((V^*)^{\otimes n} \subseteq (V^{\otimes n})^*\). However, we do not need to address the distinction between \((V^{\otimes n})^*\) and \((V^*)^{\otimes n}\), nor the subtleties concerning the dual space \(T((V))^*\), since the argument (ref.~\cref{subsec:construction-of-functionals}) for the lower bound can be carried out in a truncated tensor algebra of large enough level, while the final limit is taken solely for estimates.   
\end{remark}
\begin{example}[\protect{\cite{ChangLyonsNi2018}}]
    The common examples of reasonable tensor norms are the projective, injective and Hilbert-Schmidt tensor norms.
\end{example}

\subsection{Paths and associated Stieltjes measures}

Now, let \(\gamma \colon J=[a,b] \to V\), where \(J \subseteq \mathbb{R}\) is a compact interval, be a continuous path and let \(\mathcal{P} : a = s_0 < s_1 < ... < s_m = b \) be a partition of \(J\) into \(m+1\) points with \(m-1\) interior points. 
Then, the \emph{polygonal chain} obtained by linearly interpolating \((\gamma(s_j))_{j = 0}^m\) has the associated length
\begin{equation} \label{eq:length-of-polygonal-chain}
    L(\gamma,\mathcal{P}) := \sum_{j = 1}^m \|\gamma(s_j) - \gamma(s_{j-1})\|.
\end{equation}
The length (or variation) of the path \(\gamma\) is then defined to be the least upper bound over all partitions \(\mathcal{P}\) of \(J=[a,b]\),
\begin{equation} \label{eq:length-definition}
    L(\gamma) := \sup_{\mathcal{P} : a = s_0 < s_1 < ... < s_m = b} L(\gamma,\mathcal{P}).
\end{equation}
By the definition of the supremum, it holds that for any \(\varepsilon > 0\), there exists a partition \(\mathcal{P}_{\varepsilon}\) such that
\begin{equation} \label{eq:estimate-for-length-by-polygonal-chains}
    L(\gamma,\mathcal{P}_{\varepsilon}) > L(\gamma) - \varepsilon
\end{equation}
and the interpolating polygonal chain associated to \(\mathcal{P}_\varepsilon\) is referred to as a \emph{polygonal approximation} to \(\gamma\).
We refer to any continuous \(\gamma\) with \(L(\gamma) < \infty\) as a path of \emph{bounded variation} and denote the collection of all such paths by \(BV\). Then, \(BV\) is naturally equipped with path concatenation \(\sqcup\) and path orientation reversal \(\overleftarrow{\cdot}\). \\

Let \(\gamma \colon J=[a,b] \to V \in BV\) and recall that for any surjective, monotone and continuous map \(\varphi \colon I \to J\) between compact intervals \(I,J \subset \mathbb{R}\), \(\gamma(J) = \gamma\circ \varphi(I)\) and \(L(\gamma \circ \varphi) = L(\gamma)\) hold. Let now \(L := L(\gamma)\) be the fixed length of a path \(\gamma \colon J=[a,b] \to V \in BV\). Then, it is a general fact in metric geometry (ref.~\cite{Friz_Victoir_2010_fullbook}) that any \(BV\) path admits a \emph{parametrisation to arc length} \(\gamma = \tilde{\gamma} \circ \tilde{\varphi}\) such that
\begin{equation} \label{eq:arclength-parametrisation-function}
    \tilde{\varphi} \colon J=[a,b] \to [0,L], \qquad t \mapsto L(\gamma|_{[a,t]})
\end{equation}
and \(L(\tilde{\gamma}|_{[s,t]}) = |t-s|\). \\

Furthermore, to each path of bounded variation there exists a unique non-atomic Stieltjes (Borel) measure \(\mu_{\gamma}\) over the Borel sets of \(J\), defineable on closed intervals\footnote{This choice of definition is justified, since \(\gamma\) is here assumed to be continuous and hence has no atoms. If \(\gamma\) is a Cadlág path, then the natural definition is using half-open intervals. For the continuous case both definitions coincide.} \([s,t] \subseteq J\) by
\[
    \mu_\gamma([s,t]) = \gamma(t) - \gamma(s)
\]
and satisfying
\[
    \gamma(t) = \gamma(a) + \int_{[a,t]} \mathrm{d}\mu_\gamma
\]
for any \(t \in J\). 
To any \(V\)-valued (Stieltjes Borel) measure \(\mu\) we define its associated \emph{variation measure} to be
\begin{equation} \label{eq:variation-measure-definition}
    |\mu|(A) = \sup_{\Pi} \sum_{B \in \Pi} \|\mu(B)\|,
\end{equation}
analogously to the length of a path, where \(\Pi\) is a partition of a measureable (Borel) set \(A\) into countably many pairwise disjoint measureable (Borel) subsets. In the case of a real signed measure, the variation measure is produced by taking \(V=\mathbb{R}\) and the norm \(|\cdot|\).
Then, for any \(\gamma\colon J=[a,b] \to V \in BV\) it follows, given any \([s,t] \subseteq J\), that
\[
    L(\gamma|_{[s,t]}) = |\mu_\gamma|([s,t]) = \int_{[s,t]} \mathrm{d}|\mu_\gamma|,
\]
and in particular \(L(\gamma) = |\mu_\gamma|(J)\).
\begin{remark}
    The total variation \(\mu\mapsto |\mu|(J)\) defines a norm on the vector space of \(V\)-valued measures of finite variation on \(J\).
\end{remark}
Let \(C(J)\) be the complete space of all real valued continuous functions \(f\colon J \to \mathbb R\) over a compact interval \(J=[a,b]\) equipped with its supremum norm
\[
    \|f\|_\infty := \sup_{s \in J} |f(s)|.
\]
Since \(J\) is compact, the supremum is attained.
Then, for any function \(f \in C(J)\), we obtain the elementary estimate
\begin{equation} \label{eq:elementary-inequality-total-variation-measure}
    \left\|\int_{[a,b]} f(s) \mathrm{d}\mu_\gamma(s)\right\| \leq \|f\|_{\infty} \int_{[a,b]} \mathrm{d}|\mu_\gamma| = \|f\|_{\infty} L(\gamma).
\end{equation}
We abbreviate throughout \(\mathrm{d}\mu_\gamma = \mathrm{d}\gamma\) when convenient.
\begin{remark}
    For any \(\psi \in V^*\) and any \(\gamma \colon J=[a,b] \to V \in BV\) the associated Stieltjes measure \(\mu_\gamma\) induces the real finite, non-atomic, signed measure \(\lambda_\psi = \langle \psi, \mu_\gamma\rangle\) for which we denote its measure differential, under the dual pairing notation, by \(\langle \psi, \mathrm{d}\gamma \rangle := \mathrm{d}\lambda_\psi\) throughout.
\end{remark}

\subsection{The signature of a path}

The \emph{signature transformation}, or short \emph{signature}, \(S\colon BV \to T((V))\) of a path \(\gamma \colon J=[a,b] \to V\) of bounded variation is given by the formal power series of iterated integrals
\[
    S(\gamma) = \sum_{n = 0}^\infty \int_{[a,b]} \int_{[a,s_n]} ... \int_{[a,s_2]} \mathrm{d}\gamma(s_1) \otimes ... \otimes \mathrm{d}\gamma(s_{n-1}) \otimes \mathrm{d}\gamma(s_n).
\]

Let us recall the following properties of the signature transformation:
\begin{enumerate}[label = (\alph*)]
    \item \label{prop:signature:reparametrisation-invariance} \(S(\gamma \circ\varphi) = S(\gamma)\) for any surjective, monotone and continuous map \(\varphi \colon J \to I\) between compact intervals \(I,J \subset \mathbb{R}\), and in particular for the arc length parametrisation (see the discussion surrounding \eqref{eq:arclength-parametrisation-function}),
    \item \label{prop:signature:path-reversal-iversion} \(S(\overleftarrow{\gamma}) = S(\gamma)^{-1}\), 
    \item \label{prop:Chens-identity} \(S(\gamma|_{[s,u]}) \otimes S(\gamma|_{[u,t]}) = S(\gamma|_{[s,t]})\) for all \(s,u,t \in J\) with \(s \leq u \leq t\). Equivalently, by reparametrisation invariance \ref{prop:signature:reparametrisation-invariance} above, we can reexpress \ref{prop:Chens-identity} by \(S(\gamma \sqcup \eta) = S(\gamma) \otimes S(\eta)\). This semi-group property is referred to as \emph{Chen's identity} or \emph{Chen's relation}.
\end{enumerate}
\begin{remark}
    The properties of the signature should be compared to the length of a path. In particular, the length functional \(L\colon (BV,\sqcup) \to ([0,\infty),+)\) is parametrisation invariant and satisfies for any \(\gamma,\eta \in BV\) the identity \(L(\gamma \sqcup \eta) = L(\gamma) + L(\eta)\) as well as \(L(\overleftarrow{\gamma}) = L(\gamma)\).
\end{remark}
Set throughout \(S_n := \pi_n \circ S\) for any \(n \in \mathbb{N}_0\). Then, we may write using the definition of the simplex
\[
    S_n(\gamma) = \int_{\Delta_{n}(J)} \mathrm{d}\gamma^{\otimes n} = \int_{a < s_1 < ... < s_n < b} \mathrm{d}\gamma^{\otimes n},
\]
where \(S_0(\gamma) = 1 \in \mathbb{R}\) for any \(\gamma\colon J \to V \in BV\) and compact interval \(J=[a,b] \subset \mathbb R\). \\
In addition, let us recall the standard \emph{simplex subdivision}
\begin{equation} \label{eq:Decomposition-into-simplices}
    J^n = \bigcup_{\sigma \in \mathrm{Sym}(n)} \overline{\Delta_n(J) \cdot \sigma}.
\end{equation}
Here, the right action of \(\sigma \in \mathrm{Sym}(n)\) on coordinates \((t_1,...,t_n) \in \Delta_n(J)\) is given by \((t_1,\ldots,t_n)\cdot\sigma := (t_{\sigma(1)},\ldots,t_{\sigma(n)}).\)
Furthermore, let us denote for real signed measures \(\mu_1,...,\mu_n\) their product measure by \(\mu_1 \times ... \times \mu_n\) with abbreviation \(\mu^{\times n}\) in case \(\mu := \mu_1 = ... = \mu_n\).

Then, using the triangle inequality, the Banach inequality (cf.~\ref{def:reasonable-tensor-norms:banach} of \cref{def:reasonable-tensor-norms} and equation \eqref{eq:elementary-inequality-total-variation-measure}), the variation measure definition and the simplex decomposition, where the diagonals are of vanishing measure, we obtain the universal upper bound \cite[Section 1]{BoedihardjoGeng2023}
\begin{equation} \label{eq:Universal-Upper-Bound}
    \begin{aligned}
        \|S_n(\gamma)\|_n \leq \int_{\Delta_{n}(J)} \mathrm{d}|\mu_\gamma|^{\times n} = \frac{1}{n!} \int_{J^n} \mathrm{d}|\mu_\gamma|^{\times n} = \frac{1}{n!} L(\gamma)^n.
    \end{aligned}
\end{equation}
Here, the simplex decomposition contributes the \(n! = |\mathrm{Sym}(n)|\) factor, due to the permutation invariance of \(|\mu_\gamma|^{\times n}\).
\begin{remark}\label{rem:simplex-shuffle-subdivision}
    The inequality \eqref{eq:Universal-Upper-Bound} also admits an interpretation through the shuffle product, introduced in \cref{def:shuffle-product} later on, in terms of simple dual functions \(\phi \colon J \to V^*\) with \(\|\phi(t)\|_*=1\) for all \(t \in J\). Then, by the Banach inequality and dual description, we get
    \[
        \|S_n(\gamma)\|_n
        \leq
        \sup_{\substack{\phi:J\to V^*\ \mathrm{simple}\\ \|\phi(t)\|_*=1\ \text{for all }t\in J}} \left\langle\phi^{\otimes n},S_n(\gamma)\right\rangle
        =
        \frac{1}{n!}
        \sup_{\substack{\phi:J\to V^*\ \mathrm{simple}\\ \|\phi(t)\|_*=1\ \text{for all }t\in J}} \left(\left\langle\phi,S_1(\gamma)\right\rangle\right)^n
        =
        \frac{L(\gamma)^n}{n!}.
    \]
\end{remark}

Rearranging \eqref{eq:Universal-Upper-Bound}, we define for any \(\gamma \in BV\) the sequence
\begin{equation}\label{eq:Lambda-definition}
    \Lambda_n(\gamma) := \left(n! \|S_n(\gamma)\|_n\right)^{1/n}.
\end{equation}
It is an observation that the signature is not an injective map and one may therefore seek an understanding of the quotient \(BV/\{S(\gamma) = S(\eta)\}\), or equivalently, the space of paths on which the signature is injective. This was resolved by B.~Hambly and T.~Lyons in \cite{HamblyLyons2010} for \(V = \mathbb{R}^d\) by deriving the equivalence
\begin{equation} \label{eq:tree-like-equivalence}
    S(\gamma) = S(\eta) : \Leftrightarrow \overleftarrow{\gamma} \sqcup \eta \textrm{ is tree-like},
\end{equation}
which was then later generalised in \cite{BoedihardjoGengLyonsYang2016} to weakly geometric rough paths over arbitrary real Banach spaces under a weaker notion of tree-like paths than in \cite{HamblyLyons2010}, also showing \eqref{eq:tree-like-equivalence} to hold for real Banach space valued paths. 

\begin{definition}[Tree-like \(BV\) paths, comp.~\protect{\cite{BoedihardjoGengLyonsYang2016}}]
    A path $\gamma \colon J=[a,b] \to V \in BV$ is called \emph{tree-like} if there exists an $\mathbb{R}$-tree $\tau$, a continuous map $\phi\colon J \to \tau$ with $\phi(a) = \phi(b)$ and a map $\psi\colon \tau \to V$, such that $\gamma = \psi \circ \phi$. The path \(\gamma\) is called \emph{tree-reduced} if for any interval \([s,t] \subseteq J\) the path \(\gamma|_{[s,t]}\) is not tree-like.
\end{definition}
\begin{theorem}[comp.~\protect{\cite{HamblyLyons2010}, \cite{BoedihardjoGengLyonsYang2016}}] \label{thm:tree-reduced-path}
    To any path \(\gamma \colon J \to V\) of bounded variation there exists a unique tree-reduced \(BV\) path \(\gamma_{\mathrm{red}}\) associated to \(\gamma\) which makes the corresponding map \(t \mapsto S(\gamma_{\mathrm{red}}|_{[0,t]})\) injective. In addition, \(L(\gamma_{\mathrm{red}}) \leq L(\gamma)\).
\end{theorem}
\cref{thm:tree-reduced-path} easily follows from \cite[Corollary 1.6]{HamblyLyons2010} and \cite[Remark 4.1, Lemma 4.6]{BoedihardjoGengLyonsYang2016}.

\section{The length conjecture and its proof}
\label{sec:length-conjecture}
In studying the equivalence \eqref{eq:tree-like-equivalence}, B.~Hambly and T.~Lyons obtained 
the following result.
\begin{theorem}[comp.~\protect{\cite[Theorem 7]{HamblyLyons2010}}] \label{thm:C1-length-theorem}
    Let \(\gamma\) be a continuously differentable \(\mathbb R^d\)-valued path parametrised to arc length. Then,
    \[
        \lim_{n \to \infty} \Lambda_n(\gamma) = L(\gamma)
    \]
    for the choice of Hilbert-Schmidt tensor norms.
\end{theorem}

Recall that \(V\) denotes a real Banach space with norm \(\|\cdot\|\). \cite[Conjecture 1]{ChangLyonsNi2018}, with refinement in \cite[Corollary 1]{BoedihardjoGeng2019}, then conjectured that \cref{thm:C1-length-theorem} continues to hold if \(\mathbb R^d\) is replaced by \(V\), one assumes reasonable tensor norms (ref.~to \cref{def:reasonable-tensor-norms}) and continuous differentiability is replaced by bounded variation, provided that the path is tree-reduced.\footnote{While tree-reducedness is not explicitely written in \cite{ChangLyonsNi2018}, it is necessary for the conjecture in this form. This is also the formulation adopted in the subsequent literature \cite{BoedihardjoGeng2023}.}
The existence of the limit is then a result of \cite{ChangLyonsNi2018} and \cite{BoedihardjoGeng2019} which we state in the form we utilise.

\begin{lemma}[\protect{\cite{ChangLyonsNi2018}, \cite{BoedihardjoGeng2019}}] \label{lem:existence-of-limit}
    Let \(\gamma\colon [0,L] \to V\) be a path of bounded variation and \(\gamma_{\mathrm{red}}\) be the associated tree-reduced path by \cref{thm:tree-reduced-path}. Then, the limit \(\lim_{n \to \infty} \Lambda_n(\gamma)\) exists and it holds that
    \begin{equation} \label{eq:Existence-of-limit-of-length}
        0 \leq \tilde{L}(\gamma) := \lim_{n \to \infty} \Lambda_n(\gamma) \leq L(\gamma_{\mathrm{red}})
    \end{equation}
\end{lemma}
\begin{remark}
    The proof of \cref{lem:existence-of-limit} as in \cite{ChangLyonsNi2018} yields slightly more than the existence of the limit. Indeed, the supermultiplicativity estimate of \cite{ChangLyonsNi2018}, together with Fekete's lemma and the non-vanishing result of \cite{BoedihardjoGeng2019}, gives as stated in \cite[Corollary 1]{BoedihardjoGeng2019}
    \[
        \tilde L(\gamma) = \lim_{n\to\infty}\Lambda_n(\gamma) = \sup_{n\in\mathbb N}\Lambda_n(\gamma).
    \]
    Thus, throughout the article, the limit defining \(\tilde L(\gamma)\) may equivalently be replaced by the corresponding supremum.
\end{remark}

\subsection{Technical prerequisites to the proof}

The here presented proof draws on the following two technical results.
\begin{lemma}[Moments detect the support; comp.~\protect{\cite[Lemma 4.9]{Chen_2026}}] \label{lem:maximisation-of-moments}
    Let \(\nu\) be a non-zero compactly supported finite signed measure on \(\mathbb{R}\), then
    \[
        \limsup_{N \to \infty} \left| \int_{\mathbb{R}} x^N \mathrm{d}\nu(x)\right|^{1/N} = \max\{|x| \mid x \in \operatorname{supp} \nu \}.
    \]
\end{lemma}
\begin{remark} We note, that
    \begin{enumerate}
        \item \cref{lem:maximisation-of-moments} is similar to the standard limit of \(L^p\)-norms as \(p \to \infty\) but includes the signs coming from the powers and measure itself,
        \item the usage of the limes superior in \cref{lem:maximisation-of-moments} is necessary, since one may consider \(\nu\) to be made up of two point masses for which the limit ceases to exist.
    \end{enumerate}
\end{remark}
    We provide a proof for completeness and note that \cref{lem:maximisation-of-moments} can be derived from many different perspectives. Therefore, we do not show the arguments of \cite{Chen_2026} which concludes the result by contradiction via analytic continuation through singular sets. Instead, we conclude the above \cref{lem:maximisation-of-moments} as a consequence of the real moment generating function and its growth condition.
\begin{proof}
    Consider the real moment generating function
    \[
        F(t) = \int_{\mathbb R} e^{t x} \mathrm{d}\nu(x). 
    \]
    Then, it is a classical result\footnote{A possible justification is as follows: For \(t>0\), the upper bound follows from \(|F(t)|\leq e^{t\sup\operatorname{supp}\nu}|\nu|(\mathbb R)\). Conversely, assuming \(|F(t)|=O(e^{ct})\) for some \(c<\sup\operatorname{supp}\nu\), the cutoffs \(g_n(x)=1-\exp(-e^{n(x-a)}/n)\), \(a>c\), yield \(\nu((a,\infty))=0\) for every \(a>c\). Then, expanding \(g_n\) and applying dominated convergence, we arrive at a contradiction to the support condition. The case \(t<0\) follows by reflection. Another proof of \cref{lem:maximisation-of-moments} can be found in \cref{app:another-proof-of-moment-detection}} that we can recover the endpoints of the support of \(\nu\) via
    \[
        \begin{aligned}
            \sup \operatorname{supp}\nu &= \limsup_{t \to \infty} \frac{\log |F(t)|}{t}, \\
            -\inf \operatorname{supp}\nu &= \limsup_{t \to \infty} \frac{\log |F(-t)|}{t},
        \end{aligned}
    \]
    Consequently, defining the real exponential growth rate
    \[
        \tau(F) := \limsup \frac{1}{t} \log \max \{|F(t)|,|F(-t)|\}
    \]
    we obtain
    \[
        \tau(F) = \max_{x \in \operatorname{supp}\nu} |x|
    \]
    where the maximum is attained due to the support being compact. Since \(F\) converges absolutely and locally uniformly due to the finite support of \(\nu\), we may expand
    \[
        F(t) = \sum_{N = 0}^\infty \frac{F^{(N)}(0)}{N!} t^N
    \]
    with \(M_N = F^{(N)}(0)\) for any \(N \in \mathbb N_0\). Then, define
    \[
        M := \limsup_{n \to \infty} |M_N|^{1/N}.
    \]
    Due to the limes superior, for every \(\varepsilon > 0\) there exists a \(C_\varepsilon > 0\) such that
    \[
        |M_N| \leq C_\varepsilon (M + \varepsilon)^N, \qquad \forall N\in \mathbb{N}_0.
    \]
    Consequently, for every \(t \in \mathbb R\), we have
    \[
        |F(t)| \leq C_\varepsilon \sum_{N = 0}^\infty \frac{(M+\varepsilon)^N |t|^N}{N!} = C_\varepsilon e^{(M+\varepsilon) |t|}
    \]
    and therefore
    \[
        \tau(F) \leq M + \varepsilon
    \]
    and letting \(\varepsilon\downarrow 0\), we get
    \[
        \tau(F) \leq M.
    \]
    On the other hand, the elementary estimate
    \[
        |M_N| \leq |\nu|(\mathbb{R}) \left( \max_{x \in \operatorname{supp}\nu}  |x| \right)^N
    \]
    gives after taking the \(N\)th root and the limes superior
    \[
        M \leq \tau(F),
    \]
    which finishes the proof.
\end{proof}
\begin{remark}
    The above proof of \cref{lem:maximisation-of-moments} can be understood by the classical growth theory of entire functions in the real case. Specifically, it is a Paley-Wiener type argument.
\end{remark}
The following technical result is at the heart of the proof of the length conjecture.
\begin{proposition} \label{prop:technical-proposition}
    Let \(V\) be a real Banach space and \(\gamma\colon [0,L]\to V \in BV\) be parametrised to arc length with \(L:= L(\gamma)\) and tree-reduced. Then, for every \(\varepsilon > 0\), there exist integers \(m \geq 2, \, D\geq 0\), a non-zero compactly supported finite signed real measure \(\nu\) and, for all \(r \in \{0,...,D\}\) and \(N \in \mathbb{N}_0\), bounded linear functionals \(H_{N,r} \in (V^{\otimes (N+m-1+r)})^*\) such that the following hold:
    \begin{enumerate}[label = (\roman*)]
        \item Let
        \[
            M_N = \int_{\mathbb{R}} x^N \mathrm{d}\nu(x)
        \]
        be the \(N\)th moment of \(\nu\). Then, we have the identity
        \begin{equation}\label{eq:moment-functional-identity}
            \frac{M_N}{N!} = \sum_{r = 0}^D \langle H_{N,r}, S_{N+m-1+r}(\gamma)\rangle.
        \end{equation}
        \item The functionals \(H_{N,r}\) are of polynomial growth. In particular, they satisfy
        \begin{equation} \label{eq:growth-bound-functionals}
            \|H_{N,r}\|_{N+m-1+r,*} \leq C^{D+m-1} (1+N)^{D+m-1}.
        \end{equation}
        \item The measure \(\nu\) satisfies
        \begin{equation} \label{eq:measure-support-property}
            \operatorname{supp}\nu \cap (L-\varepsilon,\infty) \neq \varnothing
        \end{equation}
    \end{enumerate}
\end{proposition}
The necessary lemmata for \cref{prop:technical-proposition} will be supplied in \cref{sec:technical-supplementaries}. In particular, the identity \eqref{eq:moment-functional-identity} is established in \eqref{eq:moments-and-functional-relation} of \cref{subsec:construction-of-functionals}. For the proof of the length conjecture, the polynomial bound \eqref{eq:growth-bound-functionals} suffices. However, the explicit construction of \(H_{N,r}\) gives the more precise binomial estimate of \cref{lem:Growth-bound-of-functionals}, from which the polynomial bound follows as \cref{cor:polynomial-growth-functionals}. Retaining this binomial estimate also gives a more precise bound for the moments, as explained in \cref{rem:binomial-moment-estimate} after the proof of the length conjecture.


\subsection{The proof of the length conjecture}

\noindent Using \cref{lem:existence-of-limit}, we can now state the full length conjecture as a theorem.
\begin{theorem}[Length conjecture] \label{thm:Length-theorem}
    Let \(\gamma\colon [0,L] \to V\) be a path of bounded variation parametrised to arc length and tree-reduced, so that \(L:=L(\gamma)\) is its fixed length. Using reasonable tensor norms, the limit \(\lim_{n \to \infty} \left(n! \|S_n(\gamma)\|_n \right)^{1/n}\) exists and
    \[
        \lim_{n \to \infty} \left(n! \|S_n(\gamma)\|_n \right)^{1/n} = L(\gamma).
    \]
\end{theorem}
By \cref{thm:tree-reduced-path}, we obtain the following corollary.
\begin{corollary}
     Let \(J=[a,b]\) be a compact real interval. Let \(\gamma\colon J \to V\) be a path of bounded variation. Using reasonable tensor norms,
    \[
        \lim_{n \to \infty} \left(n! \|S_n(\gamma)\|_n \right)^{1/n} = L(\gamma_{\mathrm{red}})
    \]
\end{corollary}

\begin{proof}[Proof of \protect{\cref{thm:Length-theorem}}]
    W.l.o.g.~let \(\gamma\colon [0,L]\to V \in BV\) be tree-reduced and parametrised to arc length with fixed length \(L:=L(\gamma)\). Define \(\Lambda_n = (n! \|S_n(\gamma)\|_n)^{1/n}\) as in \eqref{eq:Lambda-definition} before.
    Due to \cref{lem:existence-of-limit} the limit \(\tilde{L}(\gamma) := \lim_{n \to \infty} \Lambda_n(\gamma)\) exists and we may therefore take the limes superior instead. \\

    \noindent \textit{Upper bound.} \hspace{0.3em} By \eqref{eq:Universal-Upper-Bound}, it follows that \(\limsup_{n \to \infty}\Lambda_n(\gamma) \leq L(\gamma)\). \\

    \noindent\textit{Lower bound.} \hspace{0.3em} Using \cref{prop:technical-proposition}, we rearrange \eqref{eq:moment-functional-identity} to estimate
    \begin{equation} \label{eq:moment-inequality:final-argument}
        \begin{aligned}
            |M_N| &\leq N! \sum_{r = 0}^D |\langle H_{N,r}, S_{N+m-1+r}(\gamma)\rangle| \\
            &\leq N! \sum_{r = 0}^D \|H_{N,r}\|_{N+m-1+r,*} \|S_{N+m-1+r}(\gamma)\|_{N+m-1+r} \\
            &\overset{\eqref{eq:growth-bound-functionals}}{\leq} N! C^{D + m - 1} (1 + N)^{D + m - 1} \sum_{r = 0}^D \|S_{N+m-1+r}(\gamma)\|_{N+m-1+r} \\
            &\overset{\eqref{eq:Lambda-definition}}{=} C^{D+m-1} (1+N)^{D+m-1} \sum_{r = 0}^D \frac{N!}{(N+m-1+r)!} \Lambda_{N+m-1+r}(\gamma)^{N+m-1+r}.
        \end{aligned}
    \end{equation}
    Applying the existence of the limes \(\tilde{L}(\gamma) = \lim_{n \to \infty} \Lambda_n(\gamma)\) by \cref{lem:existence-of-limit}, being the same as the limes superior, let us fix some \(\kappa > 0\).
    Then, by the definition of the limes superior there exists an \(N_0 \in \mathbb{N}\), such that for all \(n \geq N_0\)
    \begin{align*}
        \Lambda_n(\gamma) \leq \tilde{L}(\gamma) + \kappa.
    \end{align*}
    Hence, using \eqref{eq:moment-inequality:final-argument} above, yields
    \begin{equation} \label{eq:moment-inequality:final-argument:cont1}
        |M_N| \leq C^{D+m-1} (1+N)^{D+m-1} \sum_{r = 0}^D \frac{N!}{(N+m-1+r)!} (\tilde{L}(\gamma) + \kappa)^{N+m-1+r}.
    \end{equation}
    Since, by construction, \(m - 1 + r \geq 0\) it follows that
    \[
        \frac{N!}{(N+m-1+r)!} \leq 1
    \]
    and therefore
    \begin{align*}
        |M_N| \leq \underbrace{C^{D+m-1} (\tilde{L}(\gamma) + \kappa)^{m-1} \sum_{r = 0}^D (\tilde{L}(\gamma) + \kappa)^{r}}_{=: K > 0} \, \, (1+ N)^{D+m-1} (\tilde{L}(\gamma) + \kappa)^{N}
    \end{align*}
    with \(K\) independent of \(N\).
    Taking the \(N\)th root of \eqref{eq:moment-inequality:final-argument:cont1} then gives the estimate
    \begin{align} \label{eq:moment-estimate-with-limit}
        |M_N|^{1/N} \leq K^{1/N} (1 + N)^{(D+m-1)/N} (\tilde{L}(\gamma) + \kappa).
    \end{align}
    Notice that both
    \[
        \lim_{N \to \infty} K^{1/N} = 1, \quad \textrm{ and } \quad \lim_{N \to \infty} (1+N)^{k/N} = 1, \quad \forall k \in \mathbb{N}
    \]
    hold and by taking the limes superior in \(N\) of the left hand side of \eqref{eq:moment-estimate-with-limit}, we loose all dependency on other parameters, specifically the dependency on \(m\) and \(D\), and therefore the following next steps are justified.\\
    Using \eqref{eq:measure-support-property} of \cref{prop:technical-proposition} together with \cref{lem:maximisation-of-moments}, we obtain
    \begin{equation} \label{eq:sandwich-equation-for-length}
        L(\gamma) - \varepsilon < \limsup_{N \to \infty} |M_N|^{1/N} \leq \tilde{L}(\gamma) + \kappa
    \end{equation}
    for any \(\varepsilon > 0\).

    But \eqref{eq:sandwich-equation-for-length} implies as \(\varepsilon,\kappa \downarrow 0\) the lower bound
    \[
        L(\gamma) \leq \tilde{L}(\gamma),
    \]
    finishing the proof.
\end{proof} 

\begin{remark}\label{rem:binomial-moment-estimate}
    As remarked after \cref{prop:technical-proposition}, the polynomial bound for \(H_{N,r}\) is sufficient for the preceding proof. One may, however, retain the binomial coefficients and the products of coefficient norms supplied by the explicit construction in \cref{lem:Growth-bound-of-functionals}. Hence, we obtain under the same chain of inequalities \eqref{eq:moment-inequality:final-argument} and \eqref{eq:moment-inequality:final-argument:cont1} the constant
    \[
        K=\frac{(\tilde{L}(\gamma) + \kappa)^{m-1}}{(m-1)!}\prod_{l=1}^{m-1}E_l(\tilde{L}(\gamma) + \kappa), \qquad E_l(t):=\sum_{r=0}^{d_l}\frac{\|h_{l,r}\|_{r,*}t^r}{r!},
    \]
    where \(E_l\) can be understood as the elementary upper bound of the signature polynomials,
    \[
        |p_l(t)|
        \leq\sum_{r=0}^{d_l}\|h_{l,r}\|_{r,*}\|S_r(\gamma|_{[0,t]})\|_r
        \leq\sum_{r=0}^{d_l}\frac{\|h_{l,r}\|_{r,*}t^r}{r!}
        =E_l(t),
    \]
    introduced in \cref{subsec:polynomial-algebra} later on. Therefore, the growth condition of the linear functionals \(H_{N,r}\) is a consequence of the underlying algebraic structure.
\end{remark}

%
%
%
%
%
%
%

\section{Proof of the technical lemma}

\label{sec:technical-supplementaries}

This section will cover the proof of \cref{prop:technical-proposition} and will supply the necessary technical construction. W.l.o.g.~we assume any \(\gamma \colon [0,L] \to V\in BV\) to be parametrised to arc length and tree-reduced, so that \(L:=L(\gamma)\) is the fixed length. For any \(m \in \mathbb N\), let us abbreviate \(\Delta_m := \Delta_{m}([0,L])\) throughout.
\begin{definition}[Shuffle algebra] \label{def:shuffle-product}
    For \(p,q\in\mathbb{N}_0\), let \(\operatorname{Sh}(p,q)\) denote the set of all permutations
    \(\sigma\in\operatorname{Sym}(p+q)\) satisfying
    \[
        \sigma(1)<\cdots<\sigma(p),
        \qquad
        \sigma(p+1)<\cdots<\sigma(p+q).
    \]
    For \(h\in (V^{\otimes p})^*\) and \(g\in (V^{\otimes q})^*\), define their
    \emph{shuffle product}
    \[
        h\shuffle g\in (V^{\otimes(p+q)})^*
    \]
    on elementary tensors by
    \[
        (h\shuffle g)(v_1\otimes\cdots\otimes v_{p+q})
        :=
        \sum_{\sigma\in\operatorname{Sh}(p,q)}
        h(v_{\sigma(1)}\otimes\cdots\otimes v_{\sigma(p)})
        g(v_{\sigma(p+1)}\otimes\cdots\otimes v_{\sigma(p+q)}).
    \]
    Extending the shuffle product bilinearly, the graded vector space
    \(
        \oplus_{n=0}^{\infty}(V^{\otimes n})^*
    \)
    becomes a unital commutative algebra, called the \emph{shuffle algebra},
    whose unit is \(1\in (V^{\otimes 0})^*=\mathbb{R}\).
\end{definition}
\begin{remark}[Shuffle identity] \label{rem:shuffle-identity}
    If \(x\in T((V))\) is the signature of a \(BV\) path \(\gamma\), then for every
    \(h\in (V^{\otimes p})^*\) and \(g\in (V^{\otimes q})^*\),
    \[
        \langle h,S_p(\gamma)\rangle
        \langle g,S_q(\gamma)\rangle
        =
        \langle h\shuffle g,S_{p+q}(\gamma)\rangle.
    \]
    In particular, we obtain for a \(\phi \in V^*\)
    \[
        \langle \phi,S_1(\gamma)\rangle^n = \langle \phi^{\shuffle n}, S_n(\gamma)\rangle = n! \langle \phi^{\otimes n}, S_n(\gamma)\rangle
    \]
\end{remark}
The subsequent technical lemmata used in the chain of inequalities in the proof of the length conjecture are a consequence of the following standard exercise for the Hahn--Banach theorem \cite[§23.18, (HB7), p.~617]{Schechter1997} in functional analysis.
\begin{lemma} \label{lem:norm-functional-association}
    Let \(E\) be a normed vector space and let \(0\neq x\in E\). Then, there exists \(\phi \in E^*\) such that
    \[
        \|\phi\|_*=1, \qquad \langle\phi, x\rangle=\|x\|.
    \]
\end{lemma}
\begin{proof}
    Define \(\phi_0\) on \(\operatorname{span}\{x\}\) by
    \[
        \phi_0(\lambda x):=\lambda\|x\|.
    \]
    Then \(\|\phi_0\|_*=1\), and Hahn--Banach \cite[§23.18, (HB7), p.~617]{Schechter1997} extends \(\phi_0\) to \(\phi\in E^*\) without increasing its norm.
\end{proof}

\noindent We organise the proof of \cref{prop:technical-proposition} in the following way:
\begin{itemize}
    \item In \cref{subsec:localisation}, we introduce the dual point of view of the length functional of \(BV\) paths, which lets us identify a real valued function \(X\) and a product domain \(R = I_1 \times ... \times I_q \subseteq [0,L]^{q}\) of open intervals for some identifiable \(q \in \mathbb N\), on which the desired lower bound holds. We then connect natural powers of \(X\) to combinations of restricted signatures.
    \item In \cref{subsec:polynomial-algebra}, we make two observations motivating the remaining construction. Integrating over the \emph{cut points}\footnote{We call the common endpoint of two closed intervals \(I\) and \(J\) of positive length their \emph{cut point} if it is their only point of intersection, so that \(I\cap J\) is a singleton.} lets us \textquote{glue} the restricted signatures back to the full signature, but introduces contributions outside \(R\), where the lower estimate need not hold. This motivates localisation by continuous cutoff functions and their approximation by signature polynomials, using \cref{lem:separation-of-points-algebra}.
    \item In \cref{subsec:measure-construction}, we construct the measures \(\sigma_j\), whose algebraic form permits gluing and whose mass is concentrated on the intervals defining \(R\). We then define two product measures \(\rho_*\) and \(\rho\), prove an approximation property in \hyperref[claim:Measure-distance]{Claim \ref*{claim:Measure-distance}}, and establish the required support property of \(\nu=X_\#\rho\) in \hyperref[claim:measure-support-for-nu]{Claim \ref*{claim:measure-support-for-nu}}.
    \item Lastly, in \cref{subsec:construction-of-functionals}, we return to the gluing observation and give the recursive construction of the bounded linear functionals \(H_{N,r}\), establishing \eqref{eq:moment-functional-identity}. We retain the binomial coefficients in the growth estimate of \cref{lem:Growth-bound-of-functionals}, and derive the polynomial bound in \cref{cor:polynomial-growth-functionals}.
\end{itemize}

\subsection{Duality and iterated integrals}
\label{subsec:localisation}

In this subsection, we develop the dual point of view of the length \(L(\gamma)\) for some \(\gamma \in BV\). First, notice that any partition \(\mathcal{P} : 0 = \tau_0 < \tau_1 < ... < \tau_m = L\) can be identified with a point \(\tau = (\tau_1,...,\tau_{m-1}) \in \Delta_{m-1}\) and hence we can rewrite
\[
    L(\gamma) = \sup_{\substack{m\geq 2,\ \tau=(\tau_1,...,\tau_{m-1}) \in \Delta_{m-1} \\ \tau_0 = 0, \, \tau_m = L}} \sum_{j = 1}^m \|\gamma(\tau_{j}) - \gamma(\tau_{j-1})\|.
\]
Fix now \(\varepsilon > 0\) and let \(\mathcal{P} : 0 = \tau_0 < \tau_1 < ... < \tau_m = L\) be a partition, suppressing the \(\varepsilon\) dependence, with associated point \(\tau = (\tau_1,...,\tau_{m-1}) \in \Delta_{m-1}\), such that
\begin{equation} \label{eq:polygonal-approximation-half}
    L(\gamma,\mathcal{P}) = \sum_{j = 1}^m \|\gamma(\tau_{j}) - \gamma(\tau_{j-1})\|  > L(\gamma) - \frac\varepsilon2.
\end{equation}
After removing any zero increments as described in the proof of \cref{lem:localisation} below and applying \cref{lem:norm-functional-association} there exists to each vector \(\gamma(\tau_j) - \gamma(\tau_{j-1})\) an associated bounded linear functional \(\phi_j := \phi_{\gamma(\tau_j) - \gamma(\tau_{j-1})}\) satisfying
\[
    \langle \phi_j, \gamma(\tau_j) - \gamma(\tau_{j-1}) \rangle = \|\gamma(\tau_j) - \gamma(\tau_{j-1})\|, \qquad \|\phi_j\|_* = 1.
\]
Then, using the identification of a partition with a point on the corresponding standard simplex, and fixing the start and end point, we introduce the map
\[
    \begin{aligned}
        X \colon \Delta_{m-1} &\to \mathbb{R} \\
        s=(s_1,...,s_{m-1}) &\mapsto \sum_{j = 1}^m \langle \phi_j, \gamma(s_j) - \gamma(s_{j-1}) \rangle
    \end{aligned}
\]
which lets us recover
\[
    X(\tau) = L(\gamma,\mathcal{P}).
\]
For fixed length \(L:=L(\gamma)\), we can then perturb each partition point \(\tau_j\) in a small neighbourhood. This lets us define by continuity of \(X\) a domain \(R\subseteq \Delta_{m-1}\) on which the estimate
\[
    X(s) > L - \varepsilon
\]
holds.
Let us record the full construction in the following lemma.
\begin{lemma}[Perturbation of a partition] \label{lem:localisation}
    Let \(0 < \varepsilon < L\) and \(\gamma\colon [0,L] \to V\in BV\). Then, there exists an integer \(m\geq 2\), a partition \(\mathcal{P} : \tau_0 = 0 < \tau_1 < ... < \tau_m = L\) and, for each \(j \in \{1,...,m\}\), functionals \(\phi_j \in V^*\) such that 
    \begin{enumerate}[label = (\roman*)]
        \item \label{lem:localisation:existence-of-functionals} for all \(j \in \{1,...,m\}\)
        \begin{enumerate}[label = (\alph*)]
            \item \(\langle\phi_j,\gamma(\tau_j) - \gamma(\tau_{j-1})\rangle = \|\gamma(\tau_j) - \gamma(\tau_{j-1})\|\),
            \item \(\langle\phi_j,\gamma(\tau_j) - \gamma(\tau_{j-1})\rangle > 0\),
            \item \(\|\phi_j\|_{*} = 1\),
        \end{enumerate}
        with
        \[
            \sum_{j = 1}^m \phi_j(\gamma(\tau_j) - \gamma(\tau_{j-1})) > L - \frac{\varepsilon}{2},
        \]
        and
        \item \label{lem:localisation:existence-of-intervals} if \(m \geq 2\), then for each \(j \in \{1,...,m-1\}\), there exists an open interval \(I_j \Subset (0,L)\) with \(\tau_j \in I_j\) (compactly contained in \((0,L)\)), such that the intervals are naturally ordered according to
        \begin{equation} \label{eq:natural-ordering}
            \sup \overline{I_{j}} \leq \inf \overline{I}_{j+1}
        \end{equation}
        and are pairwise disjoint.
    \end{enumerate}
    Define the function \(X \colon \Delta_{m-1} \to \mathbb{R}\) for any \(s = (s_1,...,s_{m-1}) \in \Delta_{m-1}\) with \(s_0 = 0\) and \(s_m = L\) by
    \begin{equation} \label{eq:special-function-X}
        X(s) = \sum_{j = 1}^m \langle\phi_j,\gamma(s_{j}) - \gamma(s_{j-1})\rangle,
    \end{equation}
    and \(R := I_1 \times ... \times I_{m-1} \subseteq \Delta_{m-1}\). Then for every \(s \in R\)
    \[
        X(s) > L - \varepsilon.
    \] 
\end{lemma}
By the arguments above \cref{lem:localisation}, we shall only provide the missing arguments for \ref{lem:localisation:existence-of-functionals} and \ref{lem:localisation:existence-of-intervals}.
\begin{proof}
    For \ref{lem:localisation:existence-of-functionals}, it remains to ensure that all increments of the chosen partition are nonzero. Suppose that \(\gamma(\tau_j)=\gamma(\tau_{j-1})\) for some \(j\in\{1,\ldots,m\}\). If \(j<m\), we remove the partition point \(\tau_j\) and discard the corresponding functional \(\phi_j\). Since
    \[
        \gamma(\tau_{j+1})-\gamma(\tau_{j-1}) =\gamma(\tau_{j+1})-\gamma(\tau_j),
    \]
    the functional \(\phi_{j+1}\) remains norming for the resulting increment. If \(j=m\), we instead remove \(\tau_{m-1}\) and discard \(\phi_m\). In this case,
    \[
        \gamma(\tau_m)-\gamma(\tau_{m-2}) =\gamma(\tau_{m-1})-\gamma(\tau_{m-2}),
    \]
    so that \(\phi_{m-1}\) remains norming for the resulting increment.

    In either case, the polygonal length and the corresponding sum of dual pairings remain unchanged, while all retained functionals still have dual norm one. Repeating this procedure finitely many times and relabelling the remaining partition points and functionals, we obtain a partition \(\mathcal P\) for which
    \[
        \langle\phi_j,\gamma(\tau_j)-\gamma(\tau_{j-1})\rangle
        =\|\gamma(\tau_j)-\gamma(\tau_{j-1})\|>0
    \]
    for every \(j\in\{1,\ldots,m\}\), and
    \[
        \sum_{j=1}^m
        \langle\phi_j,\gamma(\tau_j)-\gamma(\tau_{j-1})\rangle
        =L(\gamma,\mathcal P)
        >L-\frac{\varepsilon}{2}.
    \]

    If only one increment remains, choose an interior point \(u\) whose image differs from both endpoints and refine the partition at \(u\). Such a point exists by continuity, since the endpoint increment is nonzero. The resulting increments are nonzero and the polygonal length does not decrease. Choosing their norming functionals by \cref{lem:norm-functional-association}, we may therefore assume \(m\geq2\).

    For \ref{lem:localisation:existence-of-intervals}, let us set \(X\colon \Delta_{m-1} \to \mathbb{R}\) with \(s = (s_1,...,s_{m-1}) \in \Delta_{m-1}\) and \(s_0 = 0\) as well as \(s_m=L\) to be
    \[
        X(s) = \sum_{j = 1}^m \langle \phi_j, \gamma(s_j) - \gamma(s_{j-1})\rangle.
    \]
    Since \(\gamma\) is continuous, each \(\phi_j\) is continuous (since they are bounded) and summations are continuous over \(\mathbb{R}\), we have that \(X\) is continuous. Therefore, there exists to each \(\tau_j\) a neighbourhood (open interval) \(I_j'\) and hence a product neighbourhood \(I_1' \times ... \times I_{m-1}' \subset \Delta_{m-1}\) such that
    \[
        - \frac{\varepsilon}{2} < X(s_1,...,s_{m-1}) - L(\gamma,\mathcal{P}) < \frac{\varepsilon}{2}, \qquad \forall s \in I_1' \times ... \times I_{m-1}'.
    \]
    One can then shrink the intervals \(I_j'\) to \(I_j \subseteq I_j'\) such that \(I_{j-1} \cap I_{j} = \varnothing\) and naturally the \(I_1,...,I_{m-1}\) are ordered according to \eqref{eq:natural-ordering}. Setting \(R = I_1 \times ... \times I_{m-1}\) proves the desired existence of the domain. \\

    Thus, by the definition of \(R\) and continuity of \(X\), we obtain
    \[
        X(s) > L(\gamma,\mathcal{P}) - \frac{\varepsilon}{2} > L - \varepsilon, \qquad \forall s \in R.
    \]
\end{proof}

Using the associated Stieltjes measure \(\mu_\gamma\) and the definition of the signature, we may therefore write
\begin{equation} \label{eq:reexpressing-functional-of-polygonal-segment-by-the-signature}
    \begin{aligned}
        \langle\phi_j,\gamma(s_{j}) - \gamma(s_{j-1})\rangle &= \langle\phi_j,\mu_\gamma([s_{j-1},s_{j}])\rangle \\
        &= \left\langle \phi_j ,\, \int_{[s_{j-1},s_{j}]} \mathrm{d}\gamma\right\rangle \\
        &= \langle \phi_j, S_1(\gamma|_{[s_{j-1},s_j]})\rangle.
    \end{aligned}
\end{equation}
Thus, it follows that
\[
    X(s) = \sum_{j = 1}^m \langle \phi_j, S_1(\gamma|_{[s_{j-1},s_j]})\rangle.
\]

Let now \(N \in \mathbb{N}\). Taking the \(N\)th power of \(X\), defining the multi-index \(\boldsymbol{j} = (j_1,...,j_m) \in \mathbb{N}_0^m\) with \(|\boldsymbol{j}| = j_1 +...+j_m\), using the multinomial formula, where we bring the \(N!\) to the other side, and utilising the shuffle identity of \cref{rem:shuffle-identity}, then yields
\begin{equation} \label{eq:Nth-power-of-X}
    \begin{aligned}
        \frac{(X(s))^N}{N!} &= \sum_{|\boldsymbol{j}| = N} \frac{1}{j_1! ... j_m!} \prod_{k = 1}^m \langle \phi_k, S_1(\gamma|_{[s_{k-1},s_k]})\rangle^{j_k} \\
        &= \sum_{|\boldsymbol{j}| = N} \frac{1}{j_1! ... j_m!} \prod_{k = 1}^m \langle \phi_k^{\shuffle j_k}, S_{j_k}(\gamma|_{[s_{k-1},s_k]})\rangle \\
        &= \sum_{|\boldsymbol{j}| = N} \prod_{k = 1}^m \langle \phi_k^{\otimes j_k}, S_{j_k}(\gamma|_{[s_{k-1},s_k]})\rangle
    \end{aligned}
\end{equation}
\begin{remark}
    Another way to conclude the identity \eqref{eq:Nth-power-of-X} is using the simplex decomposition for the induced real measures \(\langle \phi_j, \mathrm{d}\gamma\rangle\).
\end{remark}
\cref{lem:localisation} then implies for any \(N \in \mathbb N\) the pointwise estimate
\begin{equation} \label{eq:pointwise-estimate}
    X(s)^N > (L - \varepsilon)^N, \qquad \forall s \in R.
\end{equation}
The identity \eqref{eq:Nth-power-of-X} therefore relates the desired lower bound to restricted signatures. However, the right-hand side still depends on the cut points \(s_1,...,s_{m-1}\), whereas \cref{prop:technical-proposition} requires functionals of the full signature. We now explain how to remove these restrictions and why this leads to the localisation and approximation argument.

\subsection{Localisation and approximation}
\label{subsec:polynomial-algebra}

Recall the setting of \cref{subsec:localisation}, and in particular the identity \eqref{eq:Nth-power-of-X}. There are two observations motivating the measure construction in \cref{subsec:measure-construction} later on.

\paragraph{Gluing the restricted signatures.}
First, notice that consecutive restricted signatures share a \emph{cut point}, e.g.~\([s_{k-1},s_k]\cap[s_{k},s_{k+1}] = \{s_k\}\). Integrating over this point against a measure induced by the path lets us combine the corresponding simplices into a single simplex. To see this, let \(m=2\), write \(s:=s_1\), fix \(j_1+j_2=N\), and let \(\lambda_\psi:=\langle\psi,\mu_\gamma\rangle\) for some \(\psi\in V^*\). Writing out the restricted signatures, using the induced measures and applying Fubini's theorem gives
\begin{equation}\label{eq:gluing-observation}
\begin{aligned}
&\int_{[0,L]}
\left\langle\phi_1^{\otimes j_1},S_{j_1}(\gamma|_{[0,s]})\right\rangle
\left\langle\phi_2^{\otimes j_2},S_{j_2}(\gamma|_{[s,L]})\right\rangle
\,\langle \psi,\mathrm d\gamma(s)\rangle
\\
&\quad=\int_{\substack{0<t_1<\cdots<t_{j_1}<s\\s<u_1<\cdots<u_{j_2}<L}}
\left\langle\phi_1^{\otimes j_1}\otimes\psi\otimes\phi_2^{\otimes j_2},
\mathrm d\gamma^{\otimes(N+1)}\right\rangle
\\
&\quad=\left\langle\phi_1^{\otimes j_1}\otimes\psi\otimes\phi_2^{\otimes j_2},S_{N+1}(\gamma)\right\rangle,
\end{aligned}
\end{equation}
where we abbreviated the integration over a product domain in the last line.
Here, the measures induced by the continuous path have no atoms, so that the endpoints and coinciding integration variables do not contribute. The same argument applies when one of the \(j_k\) is zero, with the convention \(S_0=1\).

Thus, integrating over the cut point removes the two restrictions and produces a bounded linear functional of the full signature. 

\begin{remark}\label{rem:half-shuffle-formulation-of-gluing}
    As mentioned at the end of \cref{sec:motivation-and-historical-context}, the above observation admits a corresponding formulation using half-shuffles. Following the left half-shuffle convention of \cite{SalviDiehlLyonsPreiss2023}, specialised here for some \(n \in \mathbb N\), linear functional \(u \in (V^{\otimes n})^*\) and \(\psi\in V^*\) to be
    \[
        \langle\psi\prec u,S_{n+1}(\gamma)\rangle
        =
        \int_{[0,L]}
        \langle u,S_n(\gamma|_{[0,s]})\rangle\,
        \langle\psi,\mathrm d\gamma(s)\rangle.
    \]
    Thus, \eqref{eq:gluing-observation} can be written as
    \[
        \left\langle
            \phi_2\prec\bigl(\cdots\prec\bigl(
            \phi_2\prec(\psi\prec\phi_1^{\otimes j_1})
            \bigr)\cdots\bigr),
            S_{N+1}(\gamma)
        \right\rangle,
    \]
    with \(j_2\) occurrences of \(\phi_2\) and omitting them when \(j_2=0\).
\end{remark}

\paragraph{Localising the cut points.}
The second observation concerns the lower bound. In the case \(m=2\), \cref{lem:localisation} supplies an interval \(I:=I_1\Subset(0,L)\) on which \(X(s)>L-\varepsilon\). However, the integration in \eqref{eq:gluing-observation} runs over all of \([0,L]\), and contributions from \(s\notin I\) need not preserve this estimate. This motivates the modification
\begin{equation}\label{eq:cutoff-motivation}
    \langle\psi,\mathrm d\gamma(s)\rangle
    \quad\leadsto\quad
    w(s)\langle\psi,\mathrm d\gamma(s)\rangle,
\end{equation}
where \(w\colon[0,L]\to[0,1]\) is continuous and \(\operatorname{supp}w\subset I\). The choice of \(w\) must also ensure positive mass, since \(\lambda_\psi\) is a signed measure. We will make this choice precise in \cref{lem:partition-measures} later on.

In addition, we require that the modified integral must still be expressible through the full signature. Thus, an arbitrary continuous cutoff does not supply such an expression. We therefore approximate \(w\) by a signature polynomial of the form
\[
    w(s)\approx p(s)
    =\sum_{r=0}^{d}\left\langle h_r,S_r(\gamma|_{[0,s]})\right\rangle,
\]
where \(h_r \in (V^{\otimes r})^*\) for every \(r \in \{0,...,d\}\).
The factors introduced by \(p\) can then be combined with the restricted signatures using the shuffle identity of \cref{rem:shuffle-identity}, before integrating as in \eqref{eq:gluing-observation}. We carry this out in more detail in \cref{subsec:construction-of-functionals}.

Let us recall the standard approximation result for signatures which permits this replacement of \(w\).

\begin{lemma}[Approximation algebra] \label{lem:separation-of-points-algebra}
    Let \(\gamma \colon [0,L] \to V \in BV\) be tree-reduced (and parametrised to arc length). Define the algebra of polynomial functions
    \[
        \mathcal{A}_\gamma := \left\{ t \mapsto \sum_{j=0}^{d}\left\langle h_j,S_j(\gamma|_{[0,t]})\right\rangle \mid \, h_j \in(V^{\otimes j})^{*},\, d \in \mathbb{N}_0, \, j \in \{0,...,d\} \right\}
    \]
    with operations induced by the algebra of real continuous functions \(C([0,L])\).
    Then,
    \begin{enumerate}[label = (\roman*)]
        \item \label{lem:separation-of-points-algebra:subset} \(\mathcal{A}_\gamma \subset C([0,L])\),
        \item \label{lem:separation-of-points-algebra:unital} \(\mathcal{A}_\gamma\) is a unital real algebra,
        \item \label{lem:separation-of-points-algebra:separation} \(\mathcal{A}_\gamma\) separates the points of \([0,L],\)
    \end{enumerate}
    and consequently \(\overline{\mathcal{A}_\gamma}^{\|\cdot\|_\infty} = C([0,L])\).
\end{lemma}
\begin{proof}
    \ref{lem:separation-of-points-algebra:subset} follows by definition of \(\mathcal{A}_\gamma\). \\
    For \ref{lem:separation-of-points-algebra:unital}, note that under the algebra structure of \(C([0,L])\) due to the real numbers, we need to check that \(\mathcal{A}_\gamma\) is closed as a subalgebra. Indeed, this follows by linearity of the \(\langle h_j, \cdot \rangle\) as well as the identity of \cref{rem:shuffle-identity},
    \[
        \langle h_j, S_j(\gamma)\rangle\langle h_r', S_r(\gamma)\rangle = \langle h_j \shuffle h_r', S_{r+j}(\gamma)\rangle
    \]
    for any \(h_j \in (V^{\otimes j})^*\) and \(h'_r \in (V^{\otimes r})^*\) under the identification of the duals. Since we can pick \(1 \in (V^{\otimes 0})^*\), we also get that \(1 \in \mathcal{A}_\gamma\). \\
    To see \ref{lem:separation-of-points-algebra:separation}, define the function \(G := S(\gamma|_{[0,\cdot]})\) and let \([s,t] \subseteq [0,L]\). Suppose for the sake of contradiction that there exists \(s < t\) and \(G_s = G_t\), then by Chen's relation
    \[
        G_s^{-1} \otimes G_t = S(\gamma|_{[s,t]}) = \boldsymbol{1} \Leftrightarrow \gamma|_{[s,t]} \textrm{ is tree-like}
    \]
    yielding a contradiction, since \(\gamma\) was considered tree-reduced and thus \(G\) is injective. Hence, there exists an \(n \in \mathbb{N}\) such that \(S_n(\gamma|_{[0,s]}) \neq S_n(\gamma|_{[0,t]})\). By \cref{lem:norm-functional-association} there exists\footnote{Consider for fixed \(s<t\) the element \(x = S_n(\gamma|_{[0,t]}) - S_n(\gamma|_{[0,s]})\).} an \(h \in \left(V^{\otimes n}\right)^*\) such that
    \[
        \langle h, S_n(\gamma|_{[0,s]})\rangle \neq \langle h, S_n(\gamma|_{[0,t]})\rangle. 
    \]
    The last claim follows by the Stone--Weierstrass theorem.
\end{proof}
Consequently, to each continuous cutoff \(w\) and each \(\delta>0\), there exists a \(p\in\mathcal A_\gamma\) with \(\|p-w\|_\infty<\delta\). The remaining task is to choose these approximations so that the induced measures retain positive mass on the intervals \(I_j\) as constructed in \cref{lem:localisation} with controlled contributions on their complements.

\subsection{Construction of a special measure}
\label{subsec:measure-construction}

We now implement the localisation and approximation motivated in \cref{subsec:polynomial-algebra}. For each interval \(I_j\) defining \(R\), we seek a measure of the form \(\mathrm d\sigma_j=p_j\langle\psi_j,\mathrm d\gamma\rangle\), with \(p_j \in \mathcal{A}_\gamma\) (ref.~\cref{lem:separation-of-points-algebra}) that has positive mass on \(I_j\) and controlled behavior on its complement. In this way, \(p_j\) permits the gluing while its approximation to a cutoff controls the contributions away from \(I_j\).

To control the contributions outside of each \(I_j\), let us first recall the classical property of regularity for measures on compact Hausdorff spaces. For that, let \(\mu\) be a measure on a topological space \(E\) with Borel \(\sigma\)-algebra \(\mathcal{B}(E)\). Then, \(\mu\) is a \emph{regular Borel measure} if it is \emph{inner regular}
\[  
    {\displaystyle \mu(A)=\sup\{\mu (K)\mid K\subseteq A,K{\textrm{ compact and measurable}}\}} 
\]
for every \(A\in\mathcal B(E)\), and \emph{outer regular}
\begin{equation} \label{eq:outer-regularity}
    {\displaystyle \mu (A)=\inf\{\mu (G)\mid G\supseteq A,G{\text{ open and measurable}}\}}.
\end{equation}

\begin{theorem}[\protect{\cite[Proposition 7.2.3]{Cohn2013}}.]
    Let \(X\) be a locally compact Hausdorff space that has a countable base, and let \(\mu\) be a real positive Borel measure on \(X\) that is finite on compact sets. Then \(\mu\) is regular.
\end{theorem}
\begin{remark} \label{rem:regularity-of-variation-measure}
    Since the interval \([0,L]\) is metrizeable by \(|\cdot - \cdot|\) there is a countable basis for its topology \cite[Proposition 7.1.13]{Cohn2013}.
    Hence, the variation Stieltjes measure \(|\mu_\gamma|\) associated to a \(BV\) path \(\gamma\) is a finite regular positive Borel measure.
\end{remark}

Then, recalling the domain \(R=I_1 \times ... \times I_{m-1}\) of \cref{lem:localisation} and the signature polynomial algebra of \cref{lem:separation-of-points-algebra}, we equip each \(I_j\) with a measure constructed as follows.

\begin{lemma}[Construction of measures] \label{lem:partition-measures}
    Let \(\gamma \colon [0,L] \to V \in BV\) be parametrised to arclength with \(L = L(\gamma)\) and tree-reduced. Let \(I \Subset (0,L)\) be a nonempty open interval.\footnote{\(I\) is compactly contained in \((0,L)\).} Then, \(\gamma|_I\) is non-constant and for every \(\eta > 0\), there exists a \(\psi \in V^*\) with \(\|\psi\|_* = 1\), a \(d \in \mathbb{N}\) and a signature polynomial
    \[
        p = \sum_{j = 0}^d \langle h_j, S_j(\gamma|_{[0,\cdot]})\rangle \in \mathcal{A}_\gamma
    \]
    such that the finite signed measure defined for every Borel set \(B \subset [0,L]\) by
    \[
        \sigma(B) := \int_B p \,\langle \psi, \mathrm{d}\gamma\rangle
    \]
    satisfies \(\sigma(I) > 0\) and
    \[
        |\sigma|(I^c) < \eta \, \sigma(I),
    \]
    where \(I^c = [0,L]\backslash I\) and \(|\sigma|\) denotes the variation measure of \(\sigma\).
\end{lemma}
\begin{proof}
    Since \(\gamma\colon [0,L] \to V\) is parametrised to arc length, its variation on each subinterval of positive length equals the length of that interval. Hence, \(\gamma|_I\) is non-constant. Since \(I\) is an open interval, there exists a compact interval \([r,s] \Subset I\) with \(r < s\), such that \(\gamma(r) \neq \gamma(s).\)
    With respect to the choice of the compact interval \([r,s]\), it follows by \cref{lem:norm-functional-association} that there exists a \(\psi \in V^*\) satisfying \(\|\psi\|_*=1\) and
    \[
        \psi(\gamma(s) - \gamma(r)) = \|\gamma(s) - \gamma(r)\| > 0.
    \]
    Then, let \(\lambda := \langle \psi, \mu_\gamma \rangle\) be the associated real signed measure to the Stieltjes measure \(\mu_\gamma\).
    By definition of \(\lambda\), it is a real signed Borel measure with support on \([0,L]\) together with
    \[
        \lambda([r,s]) = \psi(\gamma(s) - \gamma(r)) > 0.
    \]
    For every Borel set \(B \subset [0,L]\), we have
    \[
        |\lambda(B)| \leq |\lambda|(B) \leq \|\psi\|_* |\mu_\gamma|([0,L]) \leq L(\gamma) < \infty
    \]
    and hence \(\lambda\) is finite.\\

    Since \(|\lambda|\) is outer regular by \cref{rem:regularity-of-variation-measure}, applied to the finite positive Borel measure \(|\lambda|\), there exists\footnote{We may take \(U:=U_\delta = (r-\delta, s+\delta)\) for some \(\delta>0\) so that \([r,s] \subset U\) and \(\overline{U_\delta}=[r-\delta,s+\delta] \subset I\) by choosing \(\delta\) accordingly, which is possible since \(I\) is an open interval.} an open set \(U\) with \([r,s] \subset U\), \(\overline{U} \subset I\) and \(|\lambda|(U\backslash [r,s]) < \theta \lambda([r,s])\) for some \(\theta \in (0,1)\). Here, the upper bound follows directly by the definition of the infimum for outer regularity \eqref{eq:outer-regularity}, since \([r,s]\) is a Borel set and hence we may choose a sufficiently small \(\tilde \delta > 0\) so that \(|\lambda|(U) < |\lambda|([r,s]) + \tilde \delta\). Rearranging yields \(|\lambda|(U\backslash [r,s]) < \tilde \delta\) and setting \(\tilde \delta = \theta \lambda([r,s])\) gives the bound we just desired. \\
    
    Since \([0,L]\) is metriseable it is a normal space \cite[Proposition 7.1.7]{Cohn2013}, \([r,s]\) and \([0,L]\backslash U\) are disjoint as well as compact, and there exists by Urysohn's lemma \cite{munkres2000section33} a continuous function \(w \colon [0,L] \to [0,1]\) such that \(w(t) = 1\) on \([r,s]\) and \(w(t)=0\) on \([0,L]\backslash U\). Hence, \(\{t \mid w(t) \neq 0\} \subseteq U\) and therefore \(\operatorname{supp} w = \overline{\{t \mid w(t)\neq 0\}} \subseteq \overline{U} \subset I\). Then,
    \begin{align}
        \int_I w \mathrm{d}\lambda &= \int_{[r,s]} w \mathrm{d}\lambda + \int_{I\backslash [r,s]} \underbrace{w}_{= 0 \textrm{ on } I\backslash U} \mathrm{d}\lambda \notag \\
        &= \lambda([r,s]) + \int_{U\backslash [r,s]} w \mathrm{d}\lambda \label{eq:expression-for-integration-over-I}
    \end{align}
    One might fear that the integral over \(U \backslash [r,s]\) might be negative. However, we can bound
    \begin{equation} \label{eq:upper-bound-of-difference-of-U-and-r-s}
        \left|\int_{U \backslash [r,s]} w \mathrm{d}\lambda \right| \leq \underbrace{\|w\|_\infty}_{\leq 1} |\lambda|(U \backslash [r,s]) < \theta \lambda([r,s])
    \end{equation}
    by the choice of \(U\). Applying the bound \eqref{eq:upper-bound-of-difference-of-U-and-r-s} to the integral over \(U\backslash [r,s]\) of \eqref{eq:expression-for-integration-over-I}, it follows that
    \begin{equation} \label{eq:sigma-integral}
        c := \int_I w\, \mathrm{d}\lambda > (1-\theta) \lambda([r,s]) > 0,
    \end{equation}
    by picking \(\theta \in (0,1)\).
    
    Since \(w\) is continuous, by \cref{lem:separation-of-points-algebra} we can approximate it using a \(p \in \mathcal{A}_\gamma\) such that \(\|p - w\|_\infty < \delta\) for some \(\delta > 0\).
    Set \(\sigma(B) := \int_B p \, \mathrm d\lambda \) for any Borel set \(B\). By construction of \(w\) and \(I^c := [0,L] \backslash I \subset [0,L]\backslash U\), we can bound
    \begin{align}
        |\sigma|(I^c) &= \int_{[0,L]\backslash I} \underbrace{|p|}_{=|p-w| \textrm{ on } [0,L]\backslash I} \mathrm{d}|\lambda| \notag \\
        &\leq \|p - w\|_\infty \, |\lambda|([0,L]\backslash I) \notag \\
        & < \delta \, |\lambda|([0,L]). \label{eq:sigma-upper-bound}
    \end{align}
    Furthermore, the definition of \(c\) in \eqref{eq:sigma-integral} and the estimate \eqref{eq:sigma-upper-bound} give
    \[
        \left|\sigma(I) - c \right| = \left|\int_I (p - w) \mathrm{d}\lambda\right| < \delta \, |\lambda|([0,L]) 
    \]
    and by choosing \(\delta > 0\) sufficiently small, it holds that
    \begin{equation} \label{eq:sigma-lower-bound-on-I}
        \sigma(I) > c - \delta \, |\lambda|([0,L]) > 0.
    \end{equation}
    Trivially, \eqref{eq:sigma-upper-bound} gives 
    \begin{equation} \label{eq:upper-bound-for-variation-sigma-on-I-complement}
        |\sigma(I^c)| \leq |\sigma|(I^c) < \delta \, |\lambda|([0,L]).
    \end{equation}
    Let now \(\eta > 0\). We want to choose the uniform approximation error \(\delta > 0\) between \(w\) and \(p\) such that
    \begin{equation} \label{eq:bound-for-polynomial-approximation-error}
        \delta \, |\lambda|([0,L]) < \eta (c - \delta \, |\lambda|([0,L])).
    \end{equation}
    Hence, \(\delta \in \left(0, \frac{c \eta}{(1+\eta) |\lambda|([0,L])}\right)\). Then, positivity of \eqref{eq:sigma-lower-bound-on-I} is preserved since
    \begin{equation} \label{eq:sigma-lower-bound-with-eta}
        \sigma(I) > c - \frac{\eta c}{(1 + \eta)} = \frac{c}{(1 + \eta)} > 0
    \end{equation}
    and combining \eqref{eq:upper-bound-for-variation-sigma-on-I-complement} with \eqref{eq:bound-for-polynomial-approximation-error} we get
    \begin{align*}
        |\sigma|(I^c) < \eta (c - \delta \, |\lambda|([0,L])) \overset{\eqref{eq:sigma-lower-bound-on-I}}{<} \eta \sigma(I)
    \end{align*}
    finishing the proof.
\end{proof}

The measure of \cref{lem:partition-measures} thus has the algebraic form required by the gluing observation \eqref{eq:gluing-observation}, together with the concentration property motivated by \eqref{eq:cutoff-motivation}. We now combine these measures over the cut points.

Fix a global error \(\eta\in(0,1)\). For each \(j\in\{1,...,m-1\}\), carry out the construction of \cref{lem:partition-measures} on \(I_j\), choosing the approximation errors jointly sufficiently small as specified in the proof of \hyperref[claim:Measure-distance]{Claim \ref*{claim:Measure-distance}} below. Set
\begin{equation}\label{eq:notation-values-of-sigma-measure}
    b_j:=\sigma_j(I_j)>0,\qquad
    B:=\prod_{j=1}^{m-1}b_j,\qquad
    a_j:=|\sigma_j|(I_j^c).
\end{equation}
Define the product measures
\begin{equation}\label{eq:definition-of-rho}
    \rho:=\chi_{\Delta_{m-1}}(\sigma_1\times\cdots\times\sigma_{m-1}),
    \qquad
    \rho_*:=(\chi_{I_1}\sigma_1)\times\cdots\times(\chi_{I_{m-1}}\sigma_{m-1})
\end{equation}
on the whole hypercube \([0,L]^{m-1}\), where \(\chi\) denotes the corresponding characteristic function. Since \(R\subsetneq\Delta_{m-1}\), the measure \(\rho_*\) is concentrated on \(R\) and \(\rho_*(R)=B>0\).

Recall that \(X>L-\varepsilon\) on \(R\). Since \(\rho_*\) is a signed measure, this pointwise estimate does not by itself give a lower bound for its moments by integration. Instead, we will use it to locate positive mass of the pushforward above \(L-\varepsilon\), and then apply \cref{lem:maximisation-of-moments}. This also permits us to replace \(\rho_*\) by the measure \(\rho\), whose form allows integration over the full simplex.

\begin{claim}\label{claim:Measure-distance}
    \[
        |\rho-\rho_*|([0,L]^{m-1})<\eta B.
    \]
\end{claim}
\begin{proof}
    Recall that \(|\rho-\rho_*|\) denotes the total variation measure of the difference. Set
    \begin{equation}\label{eq:sigma-decomposition}
        \sigma_j^*:=\chi_{I_j}\sigma_j,\qquad
        \sigma_{j,*}:=\chi_{I_j^c}\sigma_j,\qquad
        \sigma_j=\sigma_j^*+\sigma_{j,*}.
    \end{equation}
    Since \(R\subset\Delta_{m-1}\), telescoping the product gives
    \begin{equation}\label{eq:reexpression-of-distance-of-multiplicative-measures}
    \begin{aligned}
        \rho-\rho_*
        &=\chi_{\Delta_{m-1}}\left(\sigma_1\times\cdots\times\sigma_{m-1}
        -\sigma_1^*\times\cdots\times\sigma_{m-1}^*\right)\\
        &=\chi_{\Delta_{m-1}}\sum_{k=1}^{m-1}
        \sigma_1\times\cdots\times\sigma_{k-1}
        \times\sigma_{k,*}\times\sigma_{k+1}^*\times\cdots\times\sigma_{m-1}^*.
    \end{aligned}
    \end{equation}
    By the triangle inequality and multiplicativity of total variation for product measures,
    \begin{equation}\label{eq:Inequality-for-total-variation-substraction-form}
        |\rho-\rho_*|([0,L]^{m-1})
        \leq\sum_{k=1}^{m-1}|\sigma_{k,*}|(I_k^c)
        \prod_{j<k}|\sigma_j|([0,L])\prod_{j>k}|\sigma_j^*|(I_j).
    \end{equation}
    To express this bound in terms of \(B\), recall the construction of \(\sigma_j\) in the proof of \cref{lem:partition-measures} and set
    \[
        \lambda_j:=\langle\psi_j,\mu_\gamma\rangle,\qquad
        c_j:=\int_{I_j}w_j\,\mathrm d\lambda_j>0,\qquad
        \|p_j-w_j\|_\infty<\delta_j.
    \]
    Since \(0\leq w_j\leq1\), \(w_j=0\) on \(I_j^c\), and \(|\lambda_j|([0,L])\leq L\), we obtain the chain of estimates
    \begin{align*}
        |\sigma_{j,*}|(I_j^c)
        &=\int_{I_j^c}|p_j-w_j|\,\mathrm d|\lambda_j|
        \leq\delta_j|\lambda_j|([0,L]),\\
        |\sigma_j^*|(I_j)
        &\leq|\sigma_j|([0,L])
        =\int_{[0,L]}|p_j|\,\mathrm d|\lambda_j|\\
        &\leq\int_{[0,L]}(|p_j-w_j|+w_j)\,\mathrm d|\lambda_j|
        \leq(1+\delta_j)|\lambda_j|([0,L]),\\
        |\rho-\rho_*|([0,L]^{m-1})
        &\leq\sum_{k=1}^{m-1}\delta_k|\lambda_k|([0,L])
        \prod_{j\ne k}(1+\delta_j)|\lambda_j|([0,L])\\
        &\leq L^{m-1}\sum_{k=1}^{m-1}\delta_k\prod_{j\ne k}(1+\delta_j)
        \longrightarrow0
        \qquad\text{as }\max_j\delta_j\longrightarrow0.
    \end{align*}
    On the other hand,
    \[
        |b_j-c_j|
        =\left|\int_{I_j}(p_j-w_j)\,\mathrm d\lambda_j\right|
        \leq\delta_j|\lambda_j|([0,L]),
        \qquad
        B=\prod_{j=1}^{m-1}b_j\longrightarrow\prod_{j=1}^{m-1}c_j>0.
    \]
    Hence the approximation errors may be chosen jointly sufficiently small such that indeed
    \[
        |\rho-\rho_*|([0,L]^{m-1})<\eta B.
    \]
\end{proof}

Define now \(\nu:=X_\#\rho\). Then, the definition of the dual function \(X\colon \Delta_{m-1} \to \mathbb R\) found in \eqref{eq:special-function-X} extends continuously to the compact hypercube \([0,L]^{m-1}\), and in particular \(|X(s)|\leq L\) on \(\Delta_{m-1}\), since
\[
    |X(s)|\leq\sum_{j=1}^{m}\|\gamma(s_j)-\gamma(s_{j-1})\|
    \leq\sum_{j=1}^{m}(s_j-s_{j-1})=L.
\]
Consequently, \(\nu\) is a finite signed real measure with compact support in \([-L,L]\). The following claim ensures that it is nonzero and has the support required in \cref{prop:technical-proposition}.

\begin{claim}\label{claim:measure-support-for-nu}
    Let \(\eta\in(0,1)\) be the chosen global error. Then,
    \[
        \nu((L-\varepsilon,\infty))>(1-\eta)B>0.
    \]
\end{claim}
\begin{proof}
    Set \(E:=(L-\varepsilon,\infty)\). Since \(X>L-\varepsilon\) on \(R\) and \(\rho_*\) is concentrated on \(R\), we have
    \[
        (X_\#\rho_*)(E)=\rho_*(X^{-1}(E))=\rho_*(R)=B.
    \]
    Thus, by \hyperref[claim:Measure-distance]{Claim \ref*{claim:Measure-distance}},
    \begin{align*}
        |\nu(E)-B|
        &=|X_\#(\rho-\rho_*)(E)|\\
        &\leq|\rho-\rho_*|(X^{-1}(E))
        \leq|\rho-\rho_*|([0,L]^{m-1})<\eta B.
    \end{align*}
    And hence
    \[
        \nu(E) > B-\eta B=(1-\eta)B>0.
    \]
\end{proof}
Define now the moments
\begin{equation}\label{eq:moments-definition}
    M_N:=\int_{\mathbb R}x^N\,\mathrm d\nu(x)
    =\int_{\Delta_{m-1}}X(s)^N\,\mathrm d\rho(s),\qquad \forall N\in\mathbb N_0.
\end{equation}
Then, \hyperref[claim:measure-support-for-nu]{Claim \ref*{claim:measure-support-for-nu}} and \cref{lem:maximisation-of-moments} give
\[
    \operatorname{supp}\nu\cap(L-\varepsilon,\infty)\ne\varnothing,
    \qquad
    \limsup_{N\to\infty}|M_N|^{1/N}
    =\max_{x\in\operatorname{supp}\nu}|x|>L-\varepsilon.
\]
In this way, the lower estimate on \(R\) is retained through the support of \(\nu\). All that remains is to return to the gluing observation and express these moments through the full signature, with control of the resulting functionals.

\subsection{Construction of the functionals}
\label{subsec:construction-of-functionals}

Recall the gluing observation \eqref{eq:gluing-observation}. The measures constructed in \cref{subsec:measure-construction} now have the required form to carry out this procedure, while retaining the lower bound through the support of \(\nu\). We first incorporate the signature polynomials into the calculation for \(m=2\), and then repeat the argument at each cut point.

Integrating \eqref{eq:Nth-power-of-X} against \(\rho\) and using \eqref{eq:moments-definition} gives
\begin{equation}\label{eq:moments-signature-relation}
    \frac{M_N}{N!}
    =\sum_{|\boldsymbol j|=N}\int_{\Delta_{m-1}}
    \prod_{k=1}^{m}\left\langle\phi_k^{\otimes j_k},S_{j_k}(\gamma|_{[s_{k-1},s_k]})\right\rangle
    \prod_{l=1}^{m-1}\mathrm d\sigma_l(s_l).
\end{equation}
Here, by construction,
\[
    \mathrm d\sigma_l(s_l)=p_l(s_l)\langle\psi_l,\mathrm d\gamma(s_l)\rangle,
    \qquad
    p_l(s_l)=\sum_{r_l=0}^{d_l}\left\langle h_{l,r_l},S_{r_l}(\gamma|_{[0,s_l]})\right\rangle.
\]
The factors arising from \(p_l\) are coordinates of the same initial restricted signature as the factors already glued at the preceding steps. This permits the use of the shuffle identity before the next integration.

\begin{example}[A simple case]\label{ex:simple-case-functional-construction}
    Let \(m=2\), write \(s:=s_1\), and fix \(j_1+j_2=N\). Consider one term \(p(s)=\langle h_r,S_r(\gamma|_{[0,s]})\rangle\). Commuting the real factors and applying the shuffle identity, we obtain
    \begin{align*}
        J_{j_1,j_2,r} &:= \int_{[0,L]} \left\langle\phi_1^{\otimes j_1},S_{j_1}(\gamma|_{[0,s]})\right\rangle \left\langle\phi_2^{\otimes j_2},S_{j_2}(\gamma|_{[s,L]})\right\rangle \underbrace{\mathrm d\sigma(s)}_{= p(s) \langle \psi_1, \mathrm{d}\gamma(s)\rangle}\\
        &=\int_{[0,L]}
        \left\langle\phi_1^{\otimes j_1},S_{j_1}(\gamma|_{[0,s]})\right\rangle
        \left\langle\phi_2^{\otimes j_2},S_{j_2}(\gamma|_{[s,L]})\right\rangle
        \left\langle h_r,S_r(\gamma|_{[0,s]})\right\rangle
        \langle\psi_1,\mathrm d\gamma(s)\rangle\\
        &=\int_{[0,L]}
        \left\langle\phi_1^{\otimes j_1}\shuffle h_r,S_{j_1+r}(\gamma|_{[0,s]})\right\rangle
        \left\langle\phi_2^{\otimes j_2},S_{j_2}(\gamma|_{[s,L]})\right\rangle
        \langle\psi_1,\mathrm d\gamma(s)\rangle\\
        &=\int_{\substack{0<t_1<\cdots<t_{j_1+r}<s\\s<u_1<\cdots<u_{j_2}<L}}
        \left\langle(\phi_1^{\otimes j_1}\shuffle h_r)\otimes\psi_1\otimes\phi_2^{\otimes j_2},
        \mathrm d\gamma^{\otimes(N+r+1)}\right\rangle\\
        &=\left\langle(\phi_1^{\otimes j_1}\shuffle h_r)\otimes\psi_1\otimes\phi_2^{\otimes j_2},
        S_{N+r+1}(\gamma)\right\rangle.
    \end{align*}
    In the last two steps, we wrote out the restricted signatures and applied Fubini's theorem as well as combining the ordered integration variables into \(\Delta_{N+r+1}\). Hence, we have identified the bounded linear functional
    \[
        H_{j_1,j_2,r}:=(\phi_1^{\otimes j_1}\shuffle h_r)\otimes\psi_1\otimes\phi_2^{\otimes j_2},
        \qquad
        J_{j_1,j_2,r}=\langle H_{j_1,j_2,r},S_{N+r+1}(\gamma)\rangle.
    \]
\end{example}

We now apply the steps of \cref{ex:simple-case-functional-construction} inductively. For that, set \(D:=d_1+\cdots+d_{m-1}\), and introduce the multi-indices \(\boldsymbol j=(j_1,...,j_m)\) and \(\boldsymbol r=(r_1,...,r_{m-1})\), with \(|\boldsymbol j|=N\) and \(|\boldsymbol r|=r\). In the sums over \(\boldsymbol r\), we always impose \(0\leq r_l\leq d_l\) or equivalently, set \(h_{l,r_l}=0\) outside this range. After expanding \eqref{eq:moments-signature-relation}, commuting summation and integration and reindexing, we define
\begin{equation}\label{eq:integral-functional}
    J_{\boldsymbol j,\boldsymbol r}:=\int_{\Delta_{m-1}}
    \left(\prod_{k=1}^{m}\left\langle\phi_k^{\otimes j_k},S_{j_k}(\gamma|_{[s_{k-1},s_k]})\right\rangle\right)
    \prod_{l=1}^{m-1}\left\langle h_{l,r_l},S_{r_l}(\gamma|_{[0,s_l]})\right\rangle
    \langle\psi_l,\mathrm d\gamma(s_l)\rangle,
\end{equation}
so that
\begin{equation}\label{eq:moments-written-out}
\begin{aligned}
    \frac{M_N}{N!}
    &=\sum_{|\boldsymbol j|=N}\sum_{r_1=0}^{d_1}\cdots\sum_{r_{m-1}=0}^{d_{m-1}}
    J_{\boldsymbol j,\boldsymbol r}\\
    &=\sum_{r=0}^{D}\sum_{|\boldsymbol j|=N}\sum_{|\boldsymbol r|=r}J_{\boldsymbol j,\boldsymbol r}.
\end{aligned}
\end{equation}

\paragraph{Identification of functionals.}
For fixed \((\boldsymbol j,\boldsymbol r)\), define recursively, for \(t\in[0,L]\),
\begin{equation}\label{eq:initial-integral}
    J^{(1)}_{\boldsymbol j,\boldsymbol r}(t):=\int_{[0,t]}
    \left\langle\phi_1^{\otimes j_1},S_{j_1}(\gamma|_{[0,s_1]})\right\rangle
    \left\langle h_{1,r_1},S_{r_1}(\gamma|_{[0,s_1]})\right\rangle
    \left\langle\phi_2^{\otimes j_2},S_{j_2}(\gamma|_{[s_1,t]})\right\rangle
    \langle\psi_1,\mathrm d\gamma(s_1)\rangle
\end{equation}
and, for \(2\leq l\leq m-1\) where \(l\) counts the cut points \(s_1,...,s_{m-1}\),
\begin{equation}\label{eq:recursive-integrals}
    J^{(l)}_{\boldsymbol j,\boldsymbol r}(t):=\int_{[0,t]}J^{(l-1)}_{\boldsymbol j,\boldsymbol r}(s_l)
    \left\langle h_{l,r_l},S_{r_l}(\gamma|_{[0,s_l]})\right\rangle
    \left\langle\phi_{l+1}^{\otimes j_{l+1}},S_{j_{l+1}}(\gamma|_{[s_l,t]})\right\rangle
    \langle\psi_l,\mathrm d\gamma(s_l)\rangle.
\end{equation}
Fubini's theorem \cite[Section~5.2]{Cohn2013}, applied to the finite measures and bounded integrands in \eqref{eq:integral-functional}, then gives \(J_{\boldsymbol j,\boldsymbol r}=J^{(m-1)}_{\boldsymbol j,\boldsymbol r}(L)\).

The identification of the functionals follows the same calculation as in \cref{ex:simple-case-functional-construction}. For the initial step, set
\[
    H^{(1)}_{\boldsymbol j,\boldsymbol r}
    :=(\phi_1^{\otimes j_1}\shuffle h_{1,r_1})\otimes\psi_1\otimes\phi_2^{\otimes j_2},
\]
such that \eqref{eq:initial-integral} becomes
\[
    J^{(1)}_{\boldsymbol j,\boldsymbol r}(t)
    =\left\langle H^{(1)}_{\boldsymbol j,\boldsymbol r},
    S_{j_1+j_2+r_1+1}(\gamma|_{[0,t]})\right\rangle.
\]
Suppose now that the functional has been identified up to step \(l-1\), and introduce the index
\[
    q_{l-1}:=\sum_{k=1}^{l}j_k+\sum_{k=1}^{l-1}r_k+l-1.
\]
Using the shuffle identity in \eqref{eq:recursive-integrals} and the same steps as before, we obtain
\begin{align*}
    J^{(l)}_{\boldsymbol j,\boldsymbol r}(t)
    &=\int_{[0,t]}
    \left\langle H^{(l-1)}_{\boldsymbol j,\boldsymbol r}\shuffle h_{l,r_l},
    S_{q_{l-1}+r_l}(\gamma|_{[0,s_l]})\right\rangle
    \left\langle\phi_{l+1}^{\otimes j_{l+1}},S_{j_{l+1}}(\gamma|_{[s_l,t]})\right\rangle
    \langle\psi_l,\mathrm d\gamma(s_l)\rangle\\
    &=\left\langle
    (H^{(l-1)}_{\boldsymbol j,\boldsymbol r}\shuffle h_{l,r_l})\otimes\psi_l\otimes\phi_{l+1}^{\otimes j_{l+1}},
    S_{q_{l-1}+r_l+1+j_{l+1}}(\gamma|_{[0,t]})\right\rangle.
\end{align*}
Thus, the inductive step is
\begin{equation}\label{eq:inductive-step-in-construction}
    H^{(l)}_{\boldsymbol j,\boldsymbol r}
    :=(H^{(l-1)}_{\boldsymbol j,\boldsymbol r}\shuffle h_{l,r_l})\otimes\psi_l\otimes\phi_{l+1}^{\otimes j_{l+1}}.
\end{equation}
Hence, we obtain after iteration through all cut points,
\begin{equation}\label{eq:integral-functional-identification}
    J_{\boldsymbol j,\boldsymbol r}=J^{(m-1)}_{\boldsymbol j,\boldsymbol r}(L)
    =\left\langle H^{(m-1)}_{\boldsymbol j,\boldsymbol r},S_{N+m-1+r}(\gamma)\right\rangle.
\end{equation}
Plugging \eqref{eq:integral-functional-identification} into \eqref{eq:moments-written-out}, and defining
\begin{equation}\label{eq:H_N_r-definition}
    H_{N,r}:=\sum_{|\boldsymbol j|=N}\sum_{|\boldsymbol r|=r}H^{(m-1)}_{\boldsymbol j,\boldsymbol r},
\end{equation}
we obtain the identity
\begin{equation}\label{eq:moments-and-functional-relation}
\begin{aligned}
    \frac{M_N}{N!}
    &=\sum_{r=0}^{D}\sum_{|\boldsymbol j|=N}\sum_{|\boldsymbol r|=r}
    \left\langle H^{(m-1)}_{\boldsymbol j,\boldsymbol r},S_{N+m-1+r}(\gamma)\right\rangle\\
    &=\sum_{r=0}^{D}\left\langle H_{N,r},S_{N+m-1+r}(\gamma)\right\rangle.
\end{aligned}
\end{equation}
All that remains is to make sure the growth of \(H_{N,r}\) remains controlled. We first retain the binomial coefficients counting the shuffle permutations, and then derive the polynomial estimate required in \cref{prop:technical-proposition}.

\begin{lemma}\label{lem:Growth-bound-of-functionals}
    Consider the construction of \(H_{N,r}\) in \eqref{eq:H_N_r-definition}. Then, for every \(N\in\mathbb N_0\) and \(0\leq r\leq D\),
    \begin{equation}\label{eq:binomial-growth-functionals}
    \begin{aligned}
        \|H_{N,r}\|_{N+m-1+r,*}
        &\leq\binom{N+m-1}{m-1}\binom{N+m-2+r}{r}\\
        &\qquad\times\sum_{\substack{|\boldsymbol r|=r\\0\leq r_l\leq d_l}}
        \frac{r!}{r_1!\cdots r_{m-1}!}
        \prod_{l=1}^{m-1}\|h_{l,r_l}\|_{r_l,*}.
    \end{aligned}
    \end{equation}
\end{lemma}
\begin{proof}
    Fix \(\boldsymbol j=(j_1,...,j_{m})\) and \(\boldsymbol r = (r_1,...,r_{m-1})\), with \(|\boldsymbol j|=N\), \(|\boldsymbol r|=r\), and \(0\leq r_l\leq d_l\). Set
    \[
        q_l:=\sum_{k=1}^{l+1}j_k+\sum_{k=1}^{l}r_k+l,
        \qquad R_l:=\sum_{k=1}^{l}r_k,\qquad 0\leq l\leq m-1,
    \]
    where the empty sums are zero. For the initial step, use \(H^{(0)}_{\boldsymbol j,\boldsymbol r}:=\phi_1^{\otimes j_1}\), so that the same calculation applies for \(1\leq l\leq m-1\). Permutation invariance \ref{def:reasonable-tensor-norms:permutation} and the dual Banach inequality \ref{def:reasonable-tensor-norms:dual-banach} of \cref{def:reasonable-tensor-norms} give for the inductive step \eqref{eq:inductive-step-in-construction},
    \begin{align*}
        \|H^{(l)}_{\boldsymbol j,\boldsymbol r}\|_{q_l,*}
        &\leq\|H^{(l-1)}_{\boldsymbol j,\boldsymbol r}\shuffle h_{l,r_l}\|_{q_{l-1}+r_l,*}
        \|\psi_l\|_*\|\phi_{l+1}^{\otimes j_{l+1}}\|_{j_{l+1},*}\\
        &\leq\binom{q_{l-1}+r_l}{r_l}
        \|H^{(l-1)}_{\boldsymbol j,\boldsymbol r}\|_{q_{l-1},*}\|h_{l,r_l}\|_{r_l,*}.
    \end{align*}
    Here, we used \(\|\psi_l\|_*=\|\phi_k\|_*=1\) and again the dual Banach inequality to get
    \(\|\phi_k^{\otimes j_k}\|_{j_k,*}\leq\|\phi_k\|_*^{j_k}=1\).
    Since \(q_{l-1}+r_l\leq N+m-2+R_l\), iterating the estimate gives
    \begin{equation}\label{eq:uniform-estimate-H}
    \begin{aligned}
        \|H^{(m-1)}_{\boldsymbol j,\boldsymbol r}\|_{N+m-1+r,*}
        &\leq\prod_{l=1}^{m-1}
        \binom{q_{l-1}+r_l}{r_l}\|h_{l,r_l}\|_{r_l,*}\\
        &\leq\left(\prod_{l=1}^{m-1}\binom{N+m-2+R_l}{r_l}\right)
        \prod_{l=1}^{m-1}\|h_{l,r_l}\|_{r_l,*}\\
        &=\frac{(N+m-2+r)!}{(N+m-2)!\prod_{l=1}^{m-1}r_l!}
        \prod_{l=1}^{m-1}\|h_{l,r_l}\|_{r_l,*}\\
        &=\binom{N+m-2+r}{r}\frac{r!}{r_1!\cdots r_{m-1}!}
        \prod_{l=1}^{m-1}\|h_{l,r_l}\|_{r_l,*}.
    \end{aligned}
    \end{equation}
    Taking the sum in \eqref{eq:H_N_r-definition} and using the triangle inequality, we obtain
    \begin{align*}
        \|H_{N,r}\|_{N+m-1+r,*}
        &\leq\sum_{|\boldsymbol j|=N}\sum_{\substack{|\boldsymbol r|=r\\0\leq r_l\leq d_l}}
        \|H^{(m-1)}_{\boldsymbol j,\boldsymbol r}\|_{N+m-1+r,*}\\
        &\leq\binom{N+m-1}{m-1}\binom{N+m-2+r}{r}
        \sum_{\substack{|\boldsymbol r|=r\\0\leq r_l\leq d_l}}
        \frac{r!}{r_1!\cdots r_{m-1}!}
        \prod_{l=1}^{m-1}\|h_{l,r_l}\|_{r_l,*},
    \end{align*}
    since there are \(\binom{N+m-1}{m-1}\) nonnegative multi-indices \(\boldsymbol j\) of total degree \(N\).
\end{proof}

\begin{corollary}\label{cor:polynomial-growth-functionals}
    There exists a constant \(C\geq1\), independent of \(N\) and uniform in \(0\leq r\leq D\), such that
    \[
        \|H_{N,r}\|_{N+m-1+r,*}
        \leq C^{m-1+r}(1+N)^{m-1+r}
        \leq C^{D+m-1}(1+N)^{D+m-1}.
    \]
\end{corollary}
\begin{proof}
    Set \(K:=\max\{1,m-1+D\}\). The binomial coefficients in \eqref{eq:binomial-growth-functionals} satisfy
    \begin{align*}
        \binom{N+m-1}{m-1}
        &=\frac{\prod_{i=1}^{m-1}(N+i)}{(m-1)!}
        \leq K^{m-1}(1+N)^{m-1},\\
        \binom{N+m-2+r}{r}
        &=\frac{\prod_{i=1}^{r}(N+m-2+i)}{r!}
        \leq K^r(1+N)^r.
    \end{align*}
    The remaining finite sum in \eqref{eq:binomial-growth-functionals} is independent of \(N\), and there are only finitely many \(r\in\{0,...,D\}\). Hence, for some \(C_0\geq1\) independent of \(N\) and \(r\),
    \begin{align*}
        \|H_{N,r}\|_{N+m-1+r,*}
        &\leq C_0K^{m-1+r}(1+N)^{m-1+r}\\
        &\leq C^{m-1+r}(1+N)^{m-1+r}
        \leq C^{D+m-1}(1+N)^{D+m-1},
    \end{align*}
    where \(C:=C_0K\), using \(m-1+r\geq1\). The empty product for \(r=0\) is understood to be one.
\end{proof}

The identity \eqref{eq:moments-and-functional-relation}, the polynomial estimate of \cref{cor:polynomial-growth-functionals}, and the support property of \hyperref[claim:measure-support-for-nu]{Claim \ref*{claim:measure-support-for-nu}} establish the three assertions of \cref{prop:technical-proposition}, completing its proof.

\appendix
\crefalias{section}{appendix}
\crefalias{subsection}{appendix}
\crefalias{subsubsection}{appendix}
\numberwithin{equation}{section}
\numberwithin{theorem}{section}
\numberwithin{lemma}{section}
\numberwithin{example}{section}
\numberwithin{claim}{section}
\numberwithin{definition}{section}
\numberwithin{remark}{section}

\section{Another proof for the support detection}
\label{app:another-proof-of-moment-detection}
This appendix is dedicated to the proof of the support detection lemma, ref.~\cref{lem:maximisation-of-moments}, using Bernstein polynomial approximations inspired by the classical Hamburger moment problem. In this way, we can bypass any complex analysis.\\ Recall, that \(\nu\) is a compactly supported signed measure. Write
\[
    M_N = \int_{\mathbb R} x^N\,\mathrm d\nu(x),\qquad V_\nu = \max_{x\in\operatorname{supp}\nu}|x|.
\]
By the triangle inequality, \(\limsup_{N \to \infty} |M_N|^{1/n}\le V_\nu\). Suppose then for the sake of contradiction that the inequality is strict and \(V_\nu > 0\). Otherwise, the result follows immediately. Then, there exists \(v\geq 0\) such that \(v<V_\nu\) and \(C>0\) with
\[
    |M_N|\le C v^N, \qquad \forall N\in \mathbb N_0.
\]
Up to replacing \(\nu\) by its reflection under \(x\mapsto-x\) if necessary, we may assume that \(V_\nu\in\operatorname{supp}\nu\), without changing
the absolute values of its moments. Pushing \(\nu\) forward under the map \(x\mapsto (x+V_\nu)/(2 V_\nu)\), we obtain a signed measure \(\mu\) on \([0,1]\) with \(1\in\operatorname{supp}\mu\). By the binomial theorem, it follows that
\begin{equation}\label{eq:new-moments}
    \left|\int_{[0,1]} t^n\,\mathrm d\mu(t)\right| \le C s^n,\qquad s:=\frac{V_\nu+v}{2 V_\nu}<1.
\end{equation}
Then, take \(f\in C([0,1])\) with \(\operatorname{supp}f\subseteq[1-\varepsilon,1]\). Its \(n\)th approximating Bernstein polynomial is
\[
    B_nf(t)=\sum_{k=0}^n f(k/n)\binom nk t^k(1-t)^{n-k}.
\]
Since \(f(k/n)=0\) for \(k<(1-\varepsilon)n\), expanding the remaining powers of \(1-t\) and using that \(\binom nj\le(en/j)^j\), yields
\begin{equation}\label{eq:estimate-Bernstein-polynomial}
    \begin{aligned}
        \left|\int_{[0,1]} B_nf\,\mathrm d\mu\right| &\le C\|f\|_\infty \sum_{j=0}^{\lfloor\varepsilon n\rfloor} \binom nj s^{n-j}(1+s)^j \\
        &\le C\|f\|_\infty(n+1) \left[s^{1-\varepsilon}\left(\frac{e(1+s)}{\varepsilon}\right)^\varepsilon \right]^n
    \end{aligned}
\end{equation}
Taking \(\varepsilon\) sufficiently small, the base of the \(n\)th power of the right hand side in \eqref{eq:estimate-Bernstein-polynomial} is strictly less than \(1\) (see \eqref{eq:new-moments}), since it converges to \(s<1\) as \(\varepsilon\downarrow0\).
Consequently,
\[
    \int_{[0,1]} B_nf\,\mathrm d\mu\longrightarrow0.
\]
But \(B_nf\to f\) uniformly by Bernstein's approximation theorem. Therefore, \(\int_{[0,1]} f\,\mathrm d\mu=0\) for every continuous \(f\) supported near \(1\), contradicting \(1\in\operatorname{supp}\mu\) and thus proving the lemma.

\section*{Acknowledgments}

Felix Medwed gratefully acknowledges funding by the Deutsche Forschungsgemeinschaft (DFG, German Research Foundation) – CRC/TRR 388 "Rough Analysis, Stochastic Dynamics and Related Fields" – Project ID 516748464.\\
Further, the author is indebted to Sylvie Paycha and Alexander Schmeding for numerous discussions and for their thoughtful comments, which shaped and has benefitted this work greatly. The author is particularly grateful to Fabrizio Zanello for generously taking the time to listen to a detailed presentation of an early version of this work and to discuss it at length. The author also thanks Peter K.~Friz for valuable advice and the opportunity to present this work at the IMOS (TU Berlin) at the 1st October 2026.

\printbibliography[heading=bibintoc]

\end{document}